\documentclass[english]{article}
\makeatletter
\usepackage{xy}

\arxivtrue 

\ifrevise
    \renewcommand{\remove}[1]{}
    \renewcommand{\edit}[1]{{\color{blue} #1}}
\fi

\ifarxiv
    \renewcommand{\remove}[1]{}
    \renewcommand{\edit}[1]{{#1}}
\fi

\begin{document}

\title{Survey Descent: A Multipoint Generalization of Gradient Descent for Nonsmooth Optimization}\date{}
\author{X.Y. Han\thanks{ORIE, Cornell University, Ithaca, NY. \texttt{xh332@cornell.edu}} \and A.S. Lewis\thanks{ORIE, Cornell University, Ithaca, NY. \texttt{people.orie.cornell.edu/aslewis}\newline
Research supported in part by National Science Foundation Grant DMS-2006990}}

\maketitle

\begin{abstract} 
For strongly convex objectives that are smooth, the classical theory of gradient descent ensures linear convergence relative to the number of gradient evaluations.  An analogous nonsmooth theory is challenging. Even when the objective is smooth at every iterate, the corresponding local models are unstable and the number of cutting planes invoked by traditional remedies is difficult to bound, leading to convergences guarantees that are sublinear relative to the cumulative number of gradient evaluations.  We instead propose a multipoint generalization of the gradient descent iteration for local optimization.  While designed with general objectives in mind, we are motivated by a ``max-of-smooth'' model that captures the subdifferential dimension at optimality.  We prove linear convergence when the objective is itself max-of-smooth, and experiments suggest a more general phenomenon.
\end{abstract}

\section{Introduction}\label{sec:intro}

To approach our target---a fast local iteration for nonsmooth optimization---let us recall, as motivation, the classical foundations of gradient descent (GD). For some function $h$, denote the linear approximation of $h$ centered at some differentiable iterate $\tx$ by \[\ell_{\tx}^{h}(x)\equiv h(\tx) + \nabla h(\tx)^T(x-\tx),\]
where $\nabla h\left(\cdot\right)$ is the gradient of $h$. When $h$ is smooth, given a step-size $\frac{1}{L}$, the canonical GD step $\tx^{+} = \tx{-}\frac{1}{L}\nabla h\left(\tx\right)$ produces the global minimum of the following local model.\begin{equation}
h_{\tx}^\mathrm{GD}\left(x\right)=\ell_{\tx}^{h}(x)+\frac{L}{2}\left\| x-\tx\right\| _{2}^{2}.\label{eq:GD_model}
\end{equation}
The model $h_{\tx}^\mathrm{GD}$ incorporates the local linear behavior observed around $\tx$ through the gradient $\nabla h\left(\tx\right)$ as well as the prior belief that the gradients themselves are $L$-Lipschitz functions of $x$. After iterating, one then builds a new model around $\tx^{+}$ and repeats the procedure. Indeed---on convex, $L$-smooth objectives---canonical results (as described in standard texts such as \citet[Chapter 10]{beck2017}) guarantee that every GD iteration will both reduce the objective value and move closer to the true global minimizer of $h$.

When the Hessian is available, one can alternatively consider the Newton model\begin{equation}
h_{\tx}^\mathrm{Newton}\left(x\right)=\ell_{\tx}^{h}(x)+\frac{1}{2}\left(x-\tx\right)^{T}\nabla^{2}h\left(\tx\right)\left(x-\tx\right).\label{eq:Newton_model}
\end{equation}
For first-order methods, however, the Hessian $\nabla^{2}h\left(\tx\right)$ is inaccessible. This motivates popular quasi-Newton approaches such as the Broyden-Fletcher-Goldfarb-Shanno (BFGS) algorithm (as described in standard texts such as \citet[Chapter 6]{nocedal2006}) that minimizes an approximation to \eqref{eq:Newton_model} after estimating $\nabla^{2}h\left(\tx\right)$ using only gradient differences and then performs a line-search.

For nonsmooth objectives, both of the above smoothness-hypothesizing quadratic models---and their variants---are inadequate for generating theoretically-guaranteed improvement at every step: They fail to capture the discontinuity of objective gradients. Thus, some algorithms work instead with a richer class of model functions. For example, one class of multipoint methods (discussed later in Section \ref{sec:related}) adaptively refines a lower cutting-plane model of the objective through a series of ``null steps'' around each iterate until it achieves descent in a ``serious step''.  Such methods work well under reasonable conditions, finitely-many null steps always sufficing to construct a successful serious step.  However,  analyzing convergence relative to all steps (null and serious) is challenging.  Current bounds on the cumulative number of steps remain sublinear due to the difficulty of uniformly bounding the number of ``in-between'' null-steps.

\subsection{The Survey Descent Iteration}\label{sec:SD_iter}

In this work, we will propose and analyze a new local \textit{Survey Descent} iteration for nonsmooth objectives $h$.

\begin{defn} \textbf{(Survey Descent Iteration $\boldsymbol{\left\{\left(P_i^\S\right)\right\}_{i=1}^k}$)\label{def:SD_Method}} Given a \textit{survey} of $k$ points, $\S=\left\{ s_{i}\right\} _{i=1}^{k}$, at which the nonsmooth objective $h$ is differentiable and a step-control parameter $L$, for each $i=1,{\dots},k$, define the \textit{$i$-th subproblem} $\left(P_i^\S\right)$ as follows:
\begin{align}
\min_{x}\ & \left\| x-\left(s_{i}-\frac{1}{L}\nabla h\left(s_{i}\right)\right)\right\| _{2}^{2}\label{eq:SD_objective}\\
\mathrm{s.t.}\ & \ell_{s_j}^{h}(x)+\frac{L}{2}\left\| x-s_{j}\right\| _{2}^{2}\leq \ell_{s_i}^{h}(x)\ \forall\ j\neq i.\label{eq:SD_constraints}
\end{align}
We refer to the solving of all subproblems $\left\{\left(P_i^\S\right)\right\}_{i=1}^k$ as a \textit{Survey Descent iteration}.
\end{defn}

When $\left(P_i^\S\right)$ is feasible, we denote its optimal solution as $s_i^+$. When all subproblems are feasible, we say that the entire Survey Descent iteration is \textit{feasible} and call $\S^+=\left\{s_i^+\right\}_{i=1}^k$ the \textit{outputs} of the iteration; otherwise, we say the entire iteration is \textit{infeasible}. After a feasible iteration, we would update $\S{\leftarrow}\S^{+}$ and repeat the Survey Descent iteration on the updated survey.

Observe that, if $k{=}1$, Survey Descent reduces to GD since there is then only one subproblem $\left(P_{1}^{\S}\right)$ and survey point $s_1$, the constraints \eqref{eq:SD_constraints} are then empty, and the objective \eqref{eq:SD_objective} of $\left(P_{1}^{\S}\right)$ is then minimized by $s_1{-}\frac{1}{L}\nabla h(s_1)$, which is the GD-step from $s_1$. Thus, we can consider Survey Descent a generalization of GD. 

\subsubsection{Main Results}
When the function $h$ is a maximum of $k$ smooth functions, we prove local linear convergence of the Survey Descent iteration to the minimizer under reasonable conditions (Theorems \ref{thm:q_linear}-\ref{thm:Rlin}). More precisely, we prove the following local properties of Survey Descent given input surveys sufficiently close to the minimizer of $h$:
\begin{itemize}
    \item Survey Descent iterations are always feasible.
    \item Survey Descent outputs are always unique, and $h$ remains differentiable at these outputs.
    \item Surveys converge Q-linearly to $\argmin_x h(x)$ when Survey Descent is applied repeatedly.
    \item Survey function-values converge R-linearly to $\min_x h(x)$ when Survey Descent is applied repeatedly.
\end{itemize}
Here, we follow the terminology of \citet[Chapter 9]{ortega2000iterative}. A major assumption for our development and any practical implementation is the availability of a workable choice of the survey size, $k$, a question associated with the ``active structure'' of $h$ at its minimizer. We discuss this choice throughout the exposition and present some empirical heuristics in Remark \ref{rem:init}.

\subsection{Linear Convergence and Nonsmooth Objectives}

\begin{figure}[H]
\begin{centering}
\begin{tabular}{c}
\textbf{Survey Descent on Simple Max-Function Example}\tabularnewline
$h_\mathrm{max}\left(x,y\right)=\left|x-y^{2}\right|+x^{2}+2y^{2}$\tabularnewline
\end{tabular}
\par\end{centering}
\subfloat[]{\includegraphics[width=1.4in]{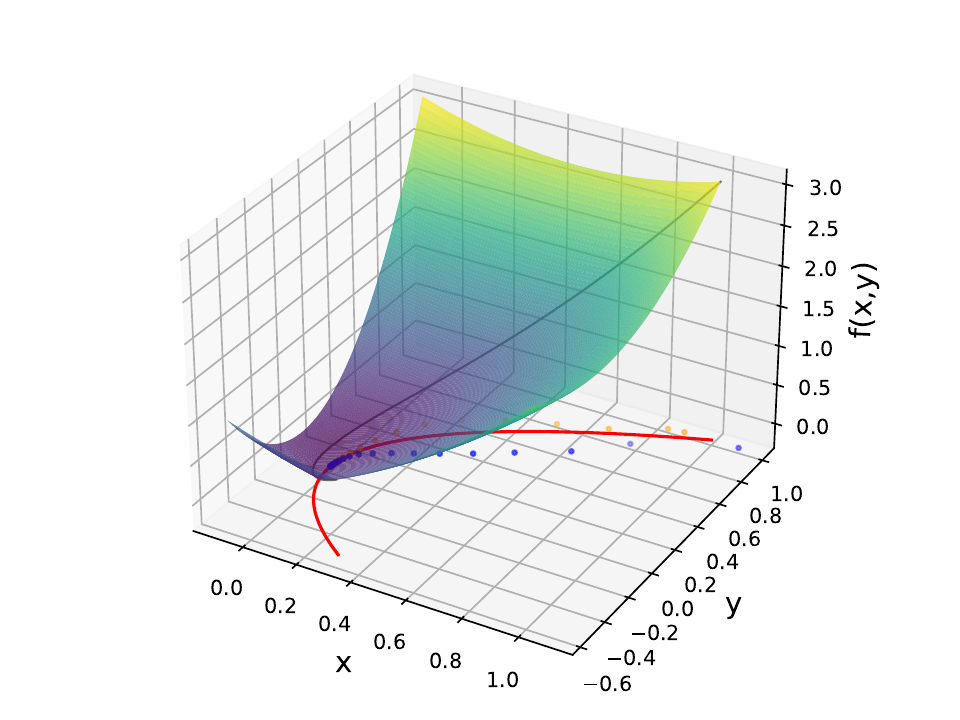}}
\subfloat[]{\includegraphics[width=1.7in]{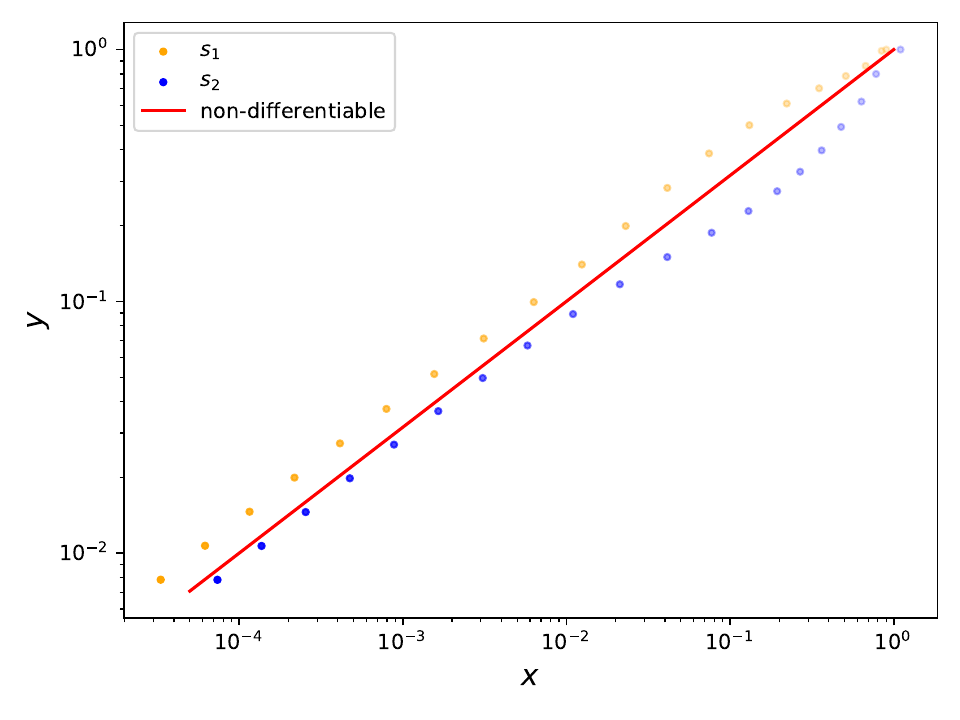}}
\subfloat[]{\includegraphics[width=1.7in]{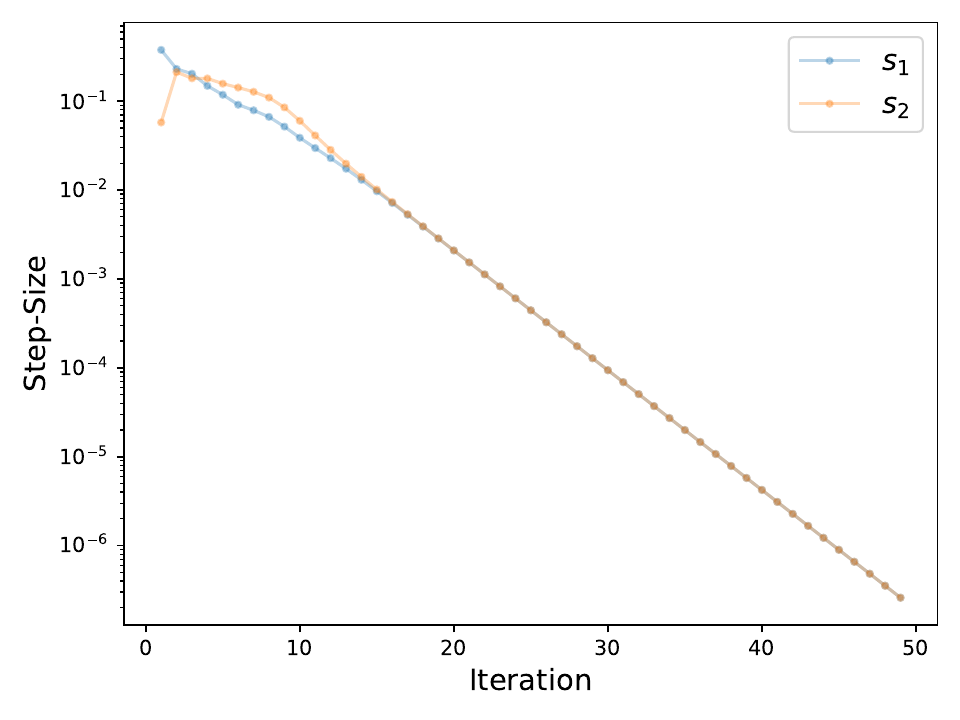}}
\subfloat[]{\includegraphics[width=1.7in]{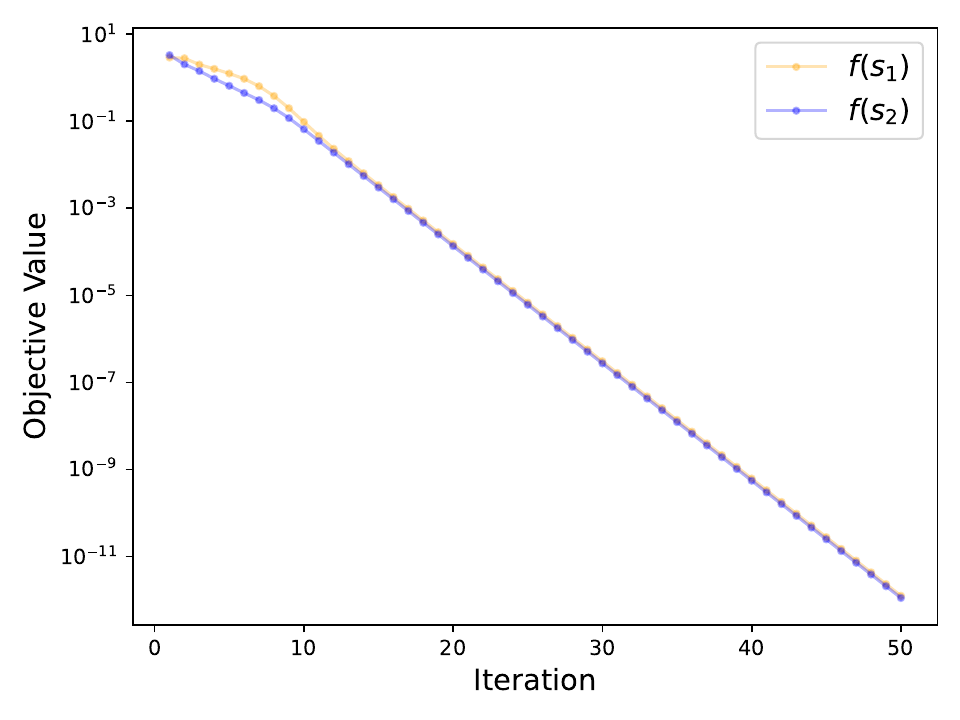}}
\caption{\textbf{\label{fig:SD_example}Survey Descent Iteration on Max-Function.} We use Survey Descent iterations to minimize the nonsmooth objective $h_\mathrm{max}\left(x,y\right)$ that achieves minimum value $0$ at unique minimizer $\bx = \left(0,0\right)$. We use step-control parameter $L = 10$ and initialize with survey points $s_{1}=\left(0.9,1\right)$ and $s_{2}=\left(1.1,1\right)$. Panel \textbf{(a)} visualizes $h_\mathrm{max}$. The Survey Descent iterates (blue and orange dots) are shown in the $xy$-plane along with the $x=y^{2}$ curve (red) on which $h_\mathrm{max}$ is non-differentiable. Panel \textbf{(b)} shows the location of iterates in the $xy$-plane where darker colors correspond to later iterations. Panel \textbf{(c)} shows the step-size (the magnitude of the difference between two consecutive iterates) at each iteration. Panel \textbf{(d)} shows the function value of the iterates. Iterates and function values converge linearly to global minimizer and minimum, respectively. Both survey points, $s_{1}$ and $s_{2}$, remain on the same smooth piece of the objective throughout the optimization. Sections \ref{sec:MOS}-\ref{sec:lconv} will theoretically derive these behaviors on objectives that are the maximum of smooth functions---of which $h_\mathrm{max}$ is a simple example.}
\end{figure}

\begin{figure}[H]
\begin{centering}
\begin{tabular}{c}
\textbf{Survey Descent in Higher Dimensions}\tabularnewline
$h_{n}\left(x\right)=\underset{{i=1,\dots,k}}{\max}\left(a_i x + x^T A_i x\right) \qquad x\in \R^n \qquad k=n/5$\tabularnewline
for random vectors $a_i$ and positive semidefinite matrices $A_i$
\end{tabular}
\par\end{centering}
\subfloat{\includegraphics[width=2.1in]{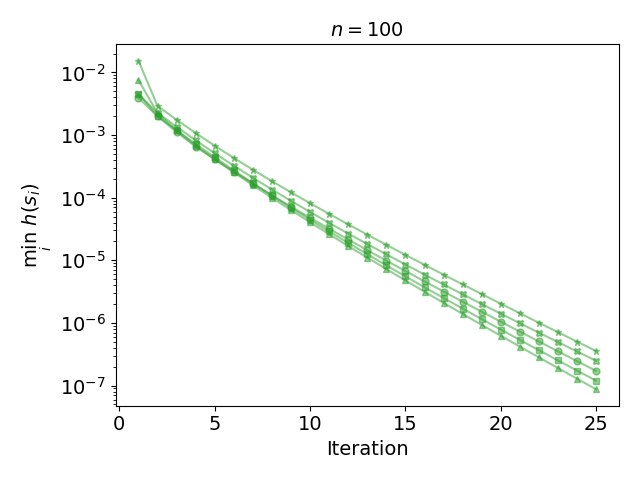}}\hfill{}\subfloat{\includegraphics[width=2.1in]{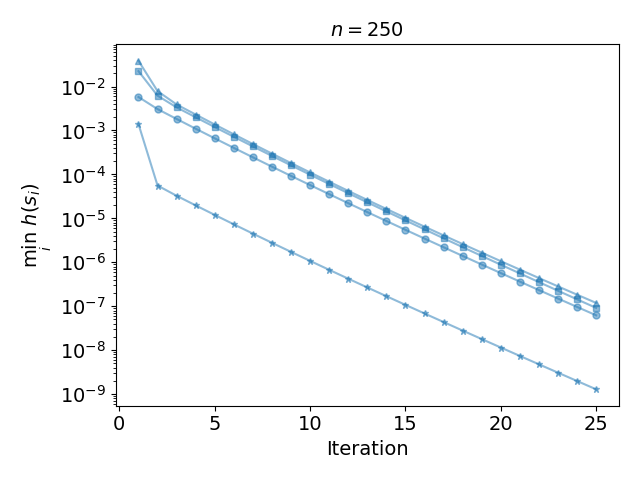}}\hfill{}\subfloat{\includegraphics[width=2.1in]{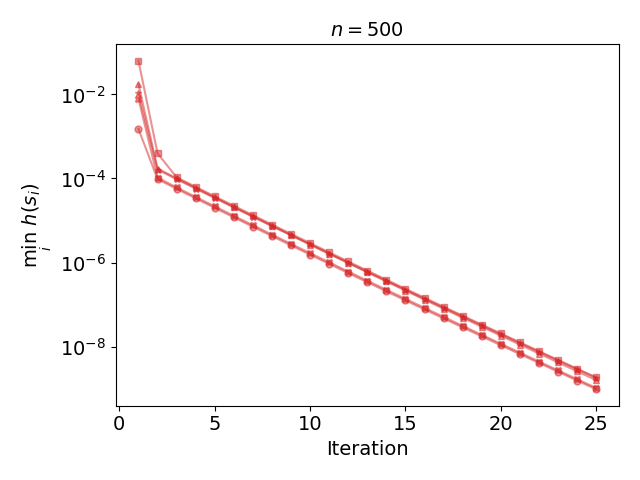}}

\caption{
\textbf{\label{fig:highdim}} Panels plot the best objective values over 25 iterations, on 5 instances, for dimension $n=100,~250,~500$. The initialization heuristic always chose the true number of components, $k$, as the survey size without \textit{a priori} information. The data $A_i = C_i^T C_i$ and $a_i$ were generated from square matrices $C_i$ and vectors with standard Gaussian entries, the set $\{a_i\}$ being adjusted to ensure that its convex hull contains zero, giving the optimal solution $x^*{=}0$.
}
\end{figure}

\begin{figure}[t]
\begin{centering}
\begin{tabular}{c}
\textbf{\textbf{Survey
Descent (with BFGS init.) on \textit{Non-}}Max-of-Smooth Objective}\tabularnewline
$h_\textrm{ME}\left(x,y,z,w\right)=\textrm{MaxEigenvalue}\left(\begin{bmatrix}x & y & z\\
y & -x & w\\
z & w & -1
\end{bmatrix}\right)$\tabularnewline
\end{tabular}
\par\end{centering}
\subfloat{\includegraphics[width=3.25in]{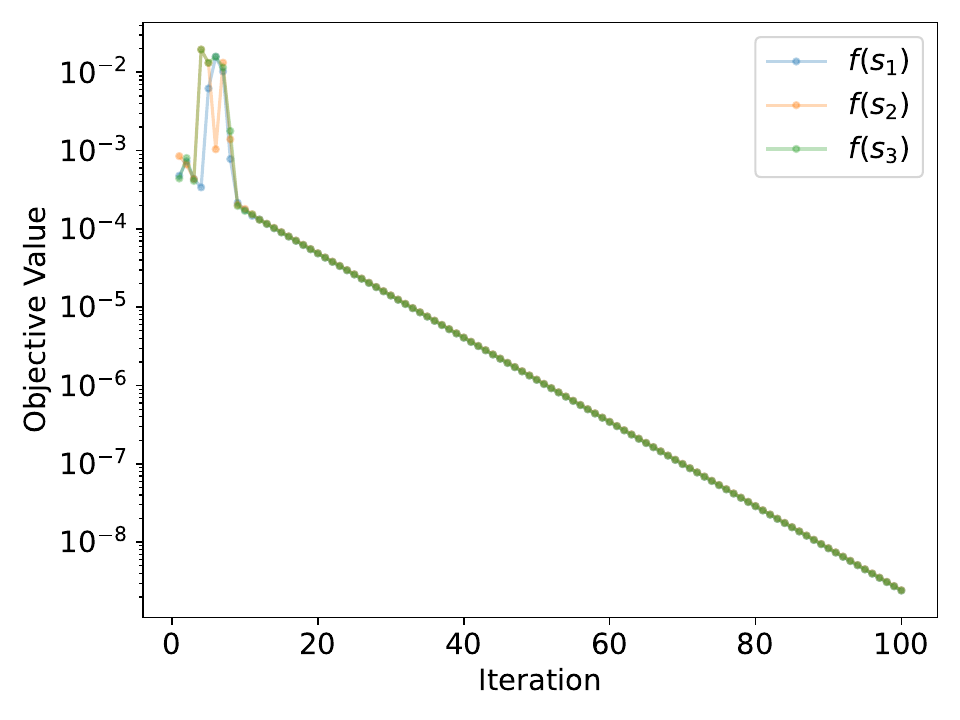}}\hfill{}\subfloat{\includegraphics[width=3.25in]{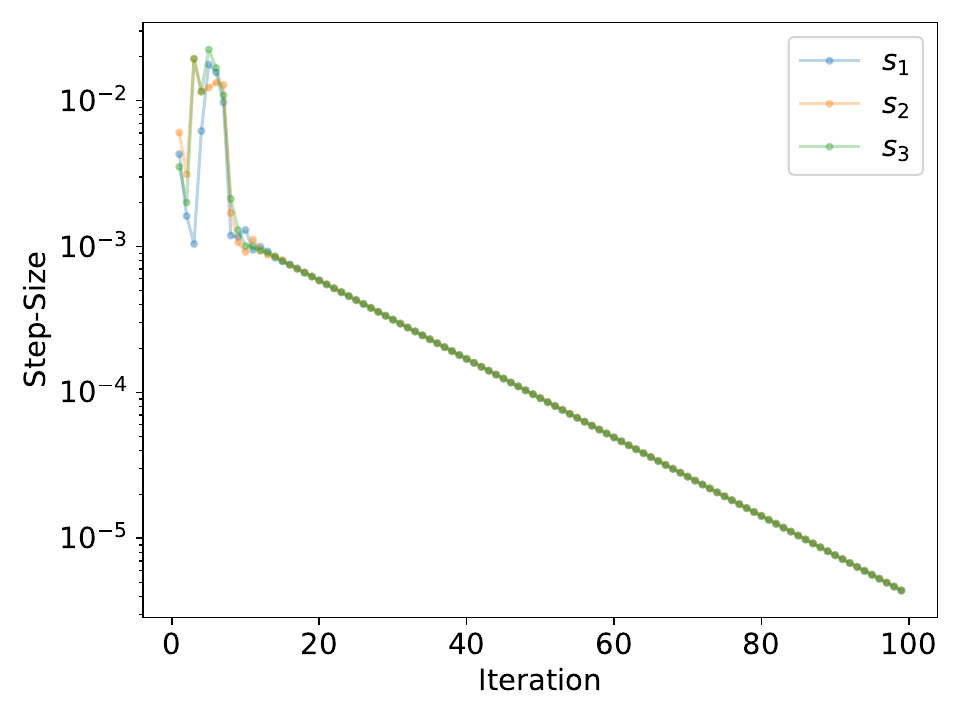}}

\caption{\textbf{\label{fig:non_max}Survey Descent of \emph{Non-}Max-of-Smooth Objective.} We use Survey Descent iterations with step-control parameter $L=10$ to minimize $h_{ME}$ with minimum value 0 at minimizer $\bx = \left(0,0,0,0\right)$. $h_{ME}$ is nonsmooth and \textit{not} expressible as a maximum of smooth functions, as is clear by noting $h_\textrm{ME}(x,y,0,0) = \sqrt{x^2+y^2}$. To create an initializing survey, we run 20 preliminary iterations of BFGS exactly as described in \citep{lewis2009} starting from $\left(1,1,1,1\right)$. We then chose both the size of a initializing survey---3 points---as well as the survey points themselves through an empirical heuristic based on the dimensionality of the BFGS iterate gradients. 
The left and right panels plot the objective values and step-sizes (the magnitude of the difference between two consecutive iterates), respectively, of each survey point during the Survey Descent iterations. After some initial fluctuations, we see the objective value and iterates of all Survey Descent survey-iterates converge at stably linear rates to the optimum. Thus, Survey Descent displays desirable performance even on an objective not captured by the theory in Sections \ref{sec:MOS}-\ref{sec:lconv}. }
\end{figure}

\begin{figure}[t]
\begin{centering}
\begin{tabular}{c}
\textbf{Survey Descent on Larger Eigenvalue Optimization Problems}\tabularnewline
$h_\textrm{ME2}\left(x\right)=\textrm{MaxEigenvalue}\left(\begin{bmatrix}0_{3\times3} & 0_{3\times7}\\
0_{7\times3} & -I_{7\times7}
\end{bmatrix}+\sum_{1}^{40}x_{i}A_{i}\right) \quad x \in \R^{40}$
\tabularnewline
for random symmetric $A_i \in \R^{10\times 10}$
\end{tabular}
\par\end{centering}
\subfloat[Objective Value]{\includegraphics[height=3.27in]{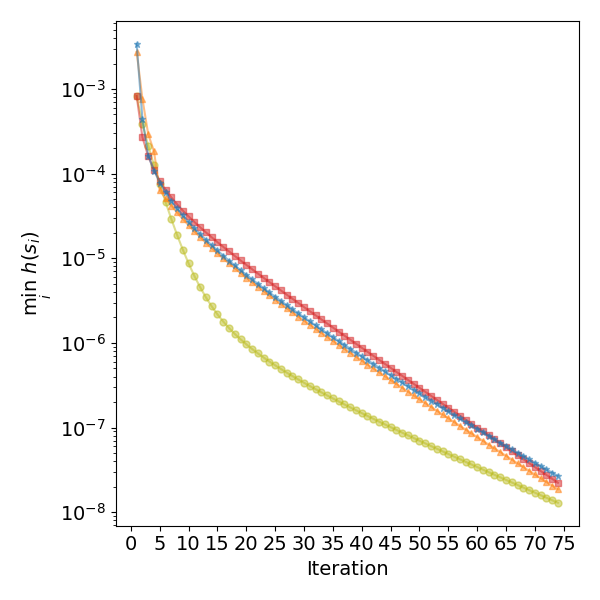}}\hfill{}\subfloat[Subdifferential Visualization]{\includegraphics[height=3.2in]{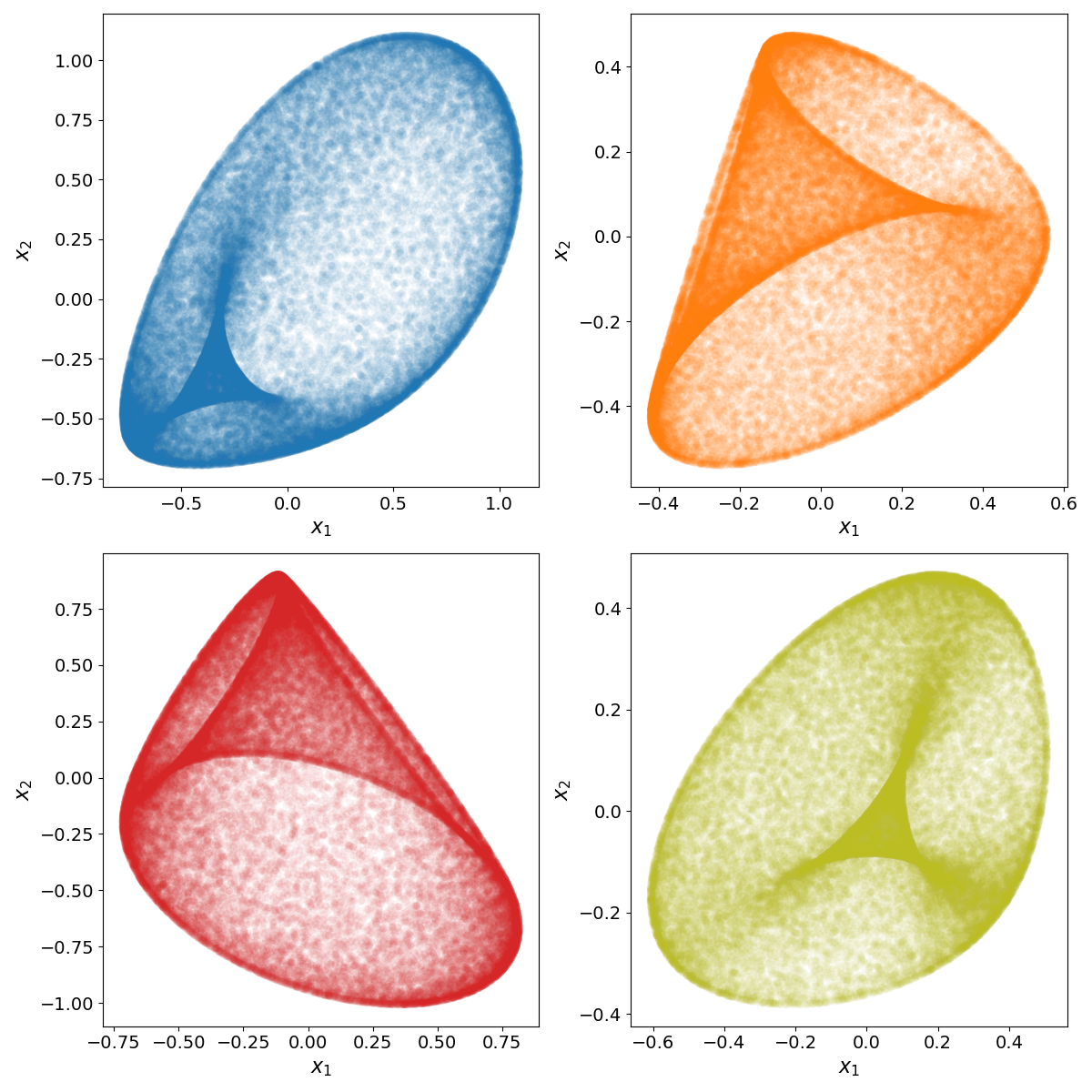}}

\caption{\textbf{\label{fig:non_max_big}Survey Descent on Larger Eigenvalue Optimization Problems.} We use Survey Descent to minimize four instances of $h_\mathrm{ME2}$. The set $\{A_i\} = \{C_i^T + C_i\}$ was generated from square matrices $C_i$ with standard Gaussian entries, shifted to ensure linear independence and an optimal solution at zero. The initialization heuristic always chose $\mathrm{dim}\left(\partial h_\mathrm{ME2}(0)\right){+}1 = 6$ as the survey size without \textit{a priori} information. Panel \textbf{(a)} plots best objective values over 75 iterations. Panel \textbf{(b)} shows two-dimensional projections of $\partial h_\mathrm{ME2}(0)$ formed by sampling $10^5$ gradients within a ball of radius $10^{-8}$. The outputs do not cluster around finitely-many points, confirming heuristically that the instances are \textit{not} max-of-smooth.}
\end{figure}

In Figure \ref{fig:SD_example}, we illustrate a simple experiment on a max-of-smooth function objective suggesting linear convergence. Figure \ref{fig:highdim} suggests that this behavior persists in higher dimensions. Moreover, Figures \ref{fig:non_max} and \ref{fig:non_max_big}, which apply Survey Descent on objectives not expressible as the maximum of smooth functions, suggest that Survey Descent may retain this linear convergence behavior even on more general nonsmooth functions, hinting at its potential as a local nonsmooth, minimization technique.

\FloatBarrier

\FloatBarrier

Many popular first-order methods exhibit empirical linear convergence on nonsmooth objectives. For example, our Figures \ref{fig:phenom} and \ref{fig:norm} illustrate the linear convergence of BFGS when minimizing two simple nonsmooth objectives, and \citet{lewis2009, lewis2013} explicitly discuss this aspect of the algorithm. Similarly, the experiments of \citet{teo2010} show that first-order multipoint methods also display linear convergence when applied to a variety of nonsmooth machine learning tasks. 

However, these behaviors are thought-provoking because modern nonsmooth, convex optimization \textit{theory} has typically only proven that canonical first-order methods converge \textit{sublinearly} when smoothness assumptions are absent---even in the presence of desirable properties such as strongly convex objectives: For example, see \citet[Chapter 3.2]{nesterov2003introductory} or \citet[Chapter 8]{beck2017} that analyzes the subgradient method as well as \citet{juditsky2011} that analyze mirror descent. Even for popular multipoint methods designed for nonsmooth optimization, previous theoretical guarantees have remained generally sublinear (discussed later in Section \ref{sec:related}). For comparison, in smooth settings, GD variants possess well-recognized \textit{linear} convergence guarantees of on $L$-smooth and $\delta$-strongly convex objectives \citep[Theorem 10.29]{beck2017}. Our derivation of local linear convergence for Survey Descent on nonsmooth objectives will directly connect to these smooth GD results.

\FloatBarrier

\subsection{Motivation of Survey Descent}

Before the formal analysis, we first describe and motivate some characteristic features within Survey Descent's design. First, note that the subproblems are simple quadratic programs with Euclidean ball constraints. Thus, they are efficiently solvable using second-order conic solvers. When all constraints moreover hold with equality, as occurs in our local convergence analysis, solving them further simplifies to routine linear algebra.

Next, the differentiability of input survey points is inspired by the common ``differentiability of every iterate'' behavior exhibited by a variety of popular first-order methods\footnote{Bundle methods (discussed in Section \ref{sec:related}) are a notable exception. Their iterates often land exactly on the nonsmooth points of their associated objectives as a consequence of minimizing their underlying cutting-plane models.} when applied to nonsmooth objectives. To illustrate, Figures \ref{fig:phenom}c and \ref{fig:norm}c show that the BFGS iterates minimizing the figures' nonsmooth objectives are always differentiable.  This behavior also manifests in \citet{lewis2009}, for example, when minimizing the nonsmooth Rosenbrock and Chebyshev-Rosenbrock objectives. Recent works have begun developing theoretical explanations for this occurrence. In particular, \citet{bianchi2020} proved that, under standard assumptions, (stochastic) gradient descent iterates are differentiable with probability-one\footnote{We distinguish the notion of ``almost all iterates of a first-order method are differentiable'' from the also common notion of ``almost all points in the domain are differentiable''---describing how all locally Lipschitz functions are non-differentiable on an at most measure-zero set in their domain. The latter result is often called Rademacher's Theorem.} even when objectives are nonconvex. Thus, we can easily imagine creating a survey of points at which $h$ is differentiable by running BFGS or GD for a few initializing iterates and collecting some subset of those iterates. (For further discussion of implementation, see Section \ref{sec:implement}.)

Lastly, the size $k$ of the survey itself relates to the dimension of the subdifferential of $h$ at its global minimizer $\bx$ (assuming it exists) while the individual points $\left\{s_{i}\right\}_{i=1}^{k}$ ``survey the landscape'' of $h$ near $\bx$. The next two subsections elaborate on these intuitions.

\FloatBarrier

\subsubsection{Dimension of the Subdifferential}\label{sec:dimsubdiff}

\begin{figure}[t]
\begin{centering}
\begin{tabular}{c}
\textbf{Simple Max-Function Example with BFGS}\tabularnewline
$h_\mathrm{max}\left(x,y\right)=\left|x-y^{2}\right|+x^{2}+2y^{2}$\tabularnewline
\end{tabular}
\par\end{centering}
\subfloat[]{\includegraphics[width=1.4in]{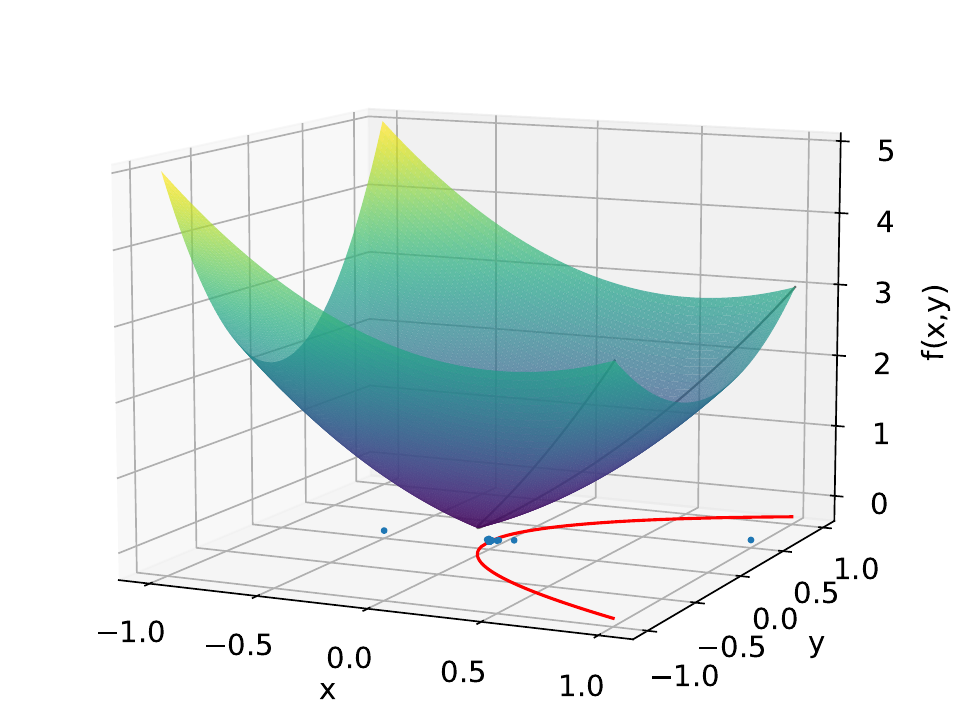}}
\subfloat[]{\includegraphics[width=1.7in]{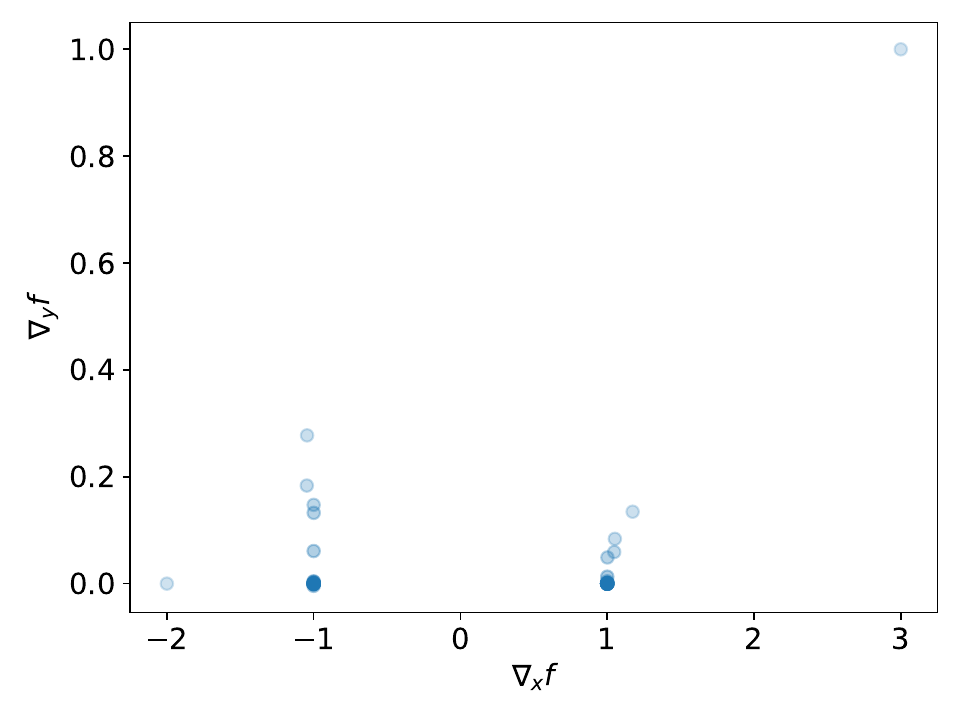}}
\subfloat[]{\includegraphics[width=1.7in]{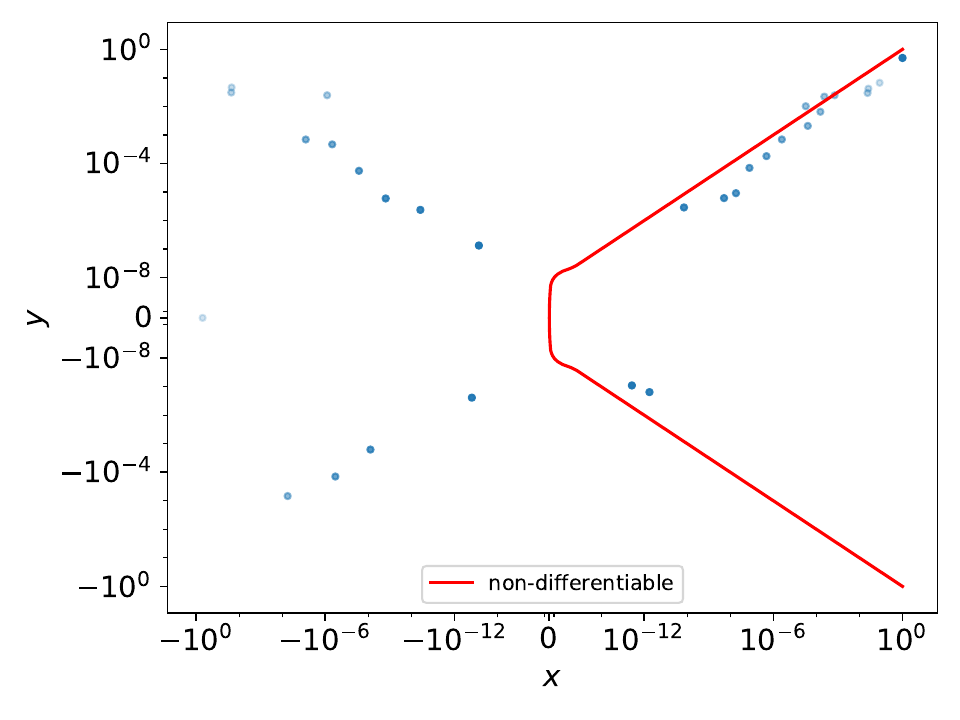}}
\subfloat[]{\includegraphics[width=1.7in]{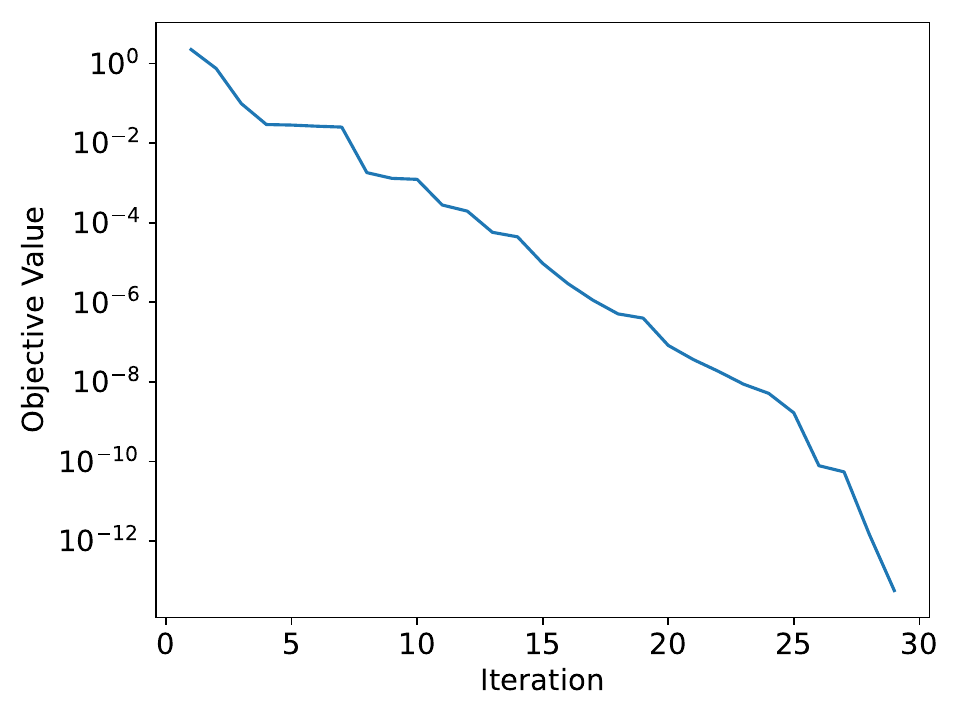}}

\caption{\textbf{\label{fig:phenom}Differentiability and Linear Convergence on Nonsmooth Objectives.} We use BFGS to minimize the nonsmooth objective $h_\mathrm{max}\left(x,y\right)$ that achieves minimum value $0$ at unique minimizer $\bx=\left(0,0\right)$. BFGS is implemented exactly as described in \citet{lewis2009}, initialized at $\left(1,0.5\right)$, and ran for 30 iterations. Panel \textbf{(a)} visualizes $h_\mathrm{max}$. It shows the iterates (blue dots) in the $xy$-plane as well as the $x=y^{2}$ curve (red) on which $h_\mathrm{max}$ is non-differentiable. Panels \textbf{(b)} and \textbf{(c)}, respectively, show the gradients at and locations of each iterate with darker colors indicating later iterations. Panel \textbf{(b)} shows gradients form a convex hull that empirically resembles $\partial h_\mathrm{max}\left(0,0\right)$'s one-dimensional structure. Panel \textbf{(c)} (in symmetric-log scale) shows that iterates lie within smooth subregions of the objective. Panel \textbf{(d)} records the linear convergence of the objective value.}
\end{figure}

The dimensionality of subdifferentials captures key nonsmoothness properties. For example, smooth functions possess everywhere zero-dimensional subdifferentials since they contain only one point: the gradient. In comparison, Figure \ref{fig:phenom} shows an objective possessing a one-dimensional subdifferential at its minimizer characterized by \begin{equation}
\partial h_\mathrm{max}\left(\bx\right)=\Conv\left(
    \left\{ \begin{bmatrix}1\\
    0
    \end{bmatrix},\begin{bmatrix}-1\\
    0
    \end{bmatrix}\right\}
\right),\label{eq:max_subdiff}
\end{equation}
where $\Conv(\cdot)$ denotes the convex hull. The objective is actually ``partly smooth'': It has a ``smooth edge'' in the $\left[0,1\right]^{T}$-direction orthogonal to $\partial h_\mathrm{max}\left(\bx\right)$. Partial smoothness is formalized and explored in detail by \citet{lewis2002} and follow-up works. 

\begin{figure}[b!]
\begin{centering}
\begin{tabular}{c}
\textbf{Simple Elliptical-Norm Example with BFGS}\tabularnewline
$h_\mathrm{ellipse}\left(x,y\right)=\sqrt{x^{2}+2y^{2}}$\tabularnewline
\end{tabular}
\par\end{centering}
\subfloat[]{\includegraphics[width=1.4in]{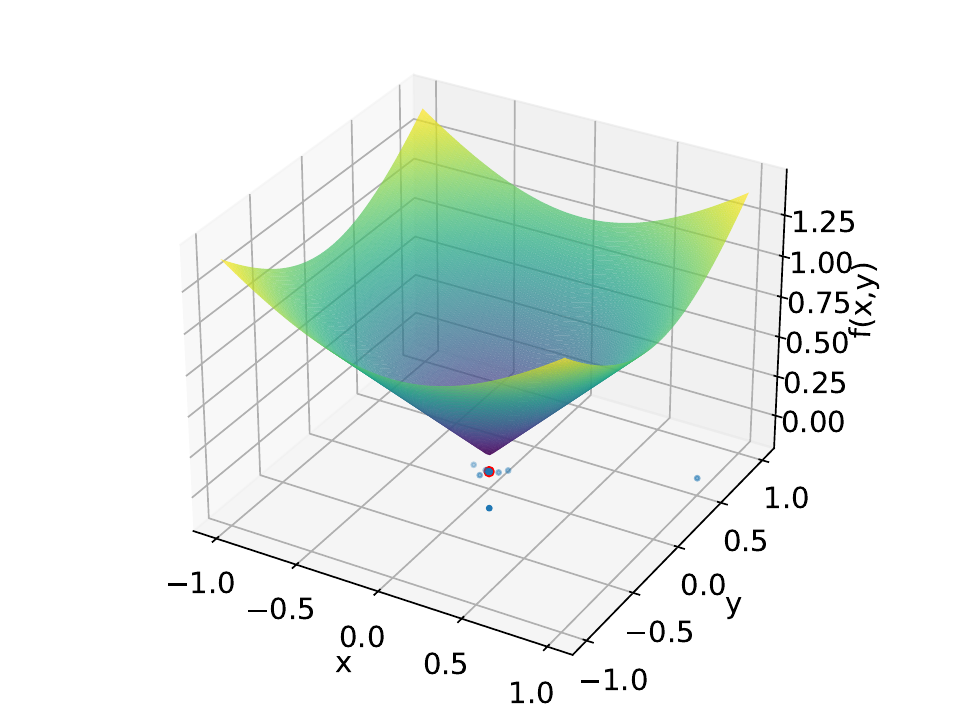}}
\subfloat[]{\includegraphics[width=1.7in]{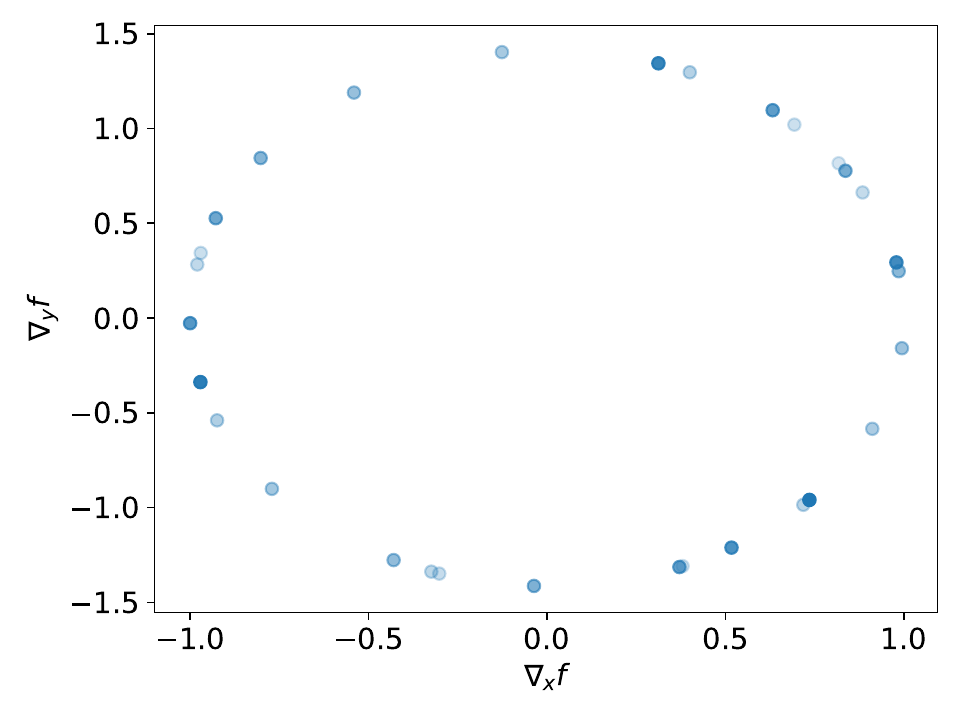}}
\subfloat[]{\includegraphics[width=1.7in]{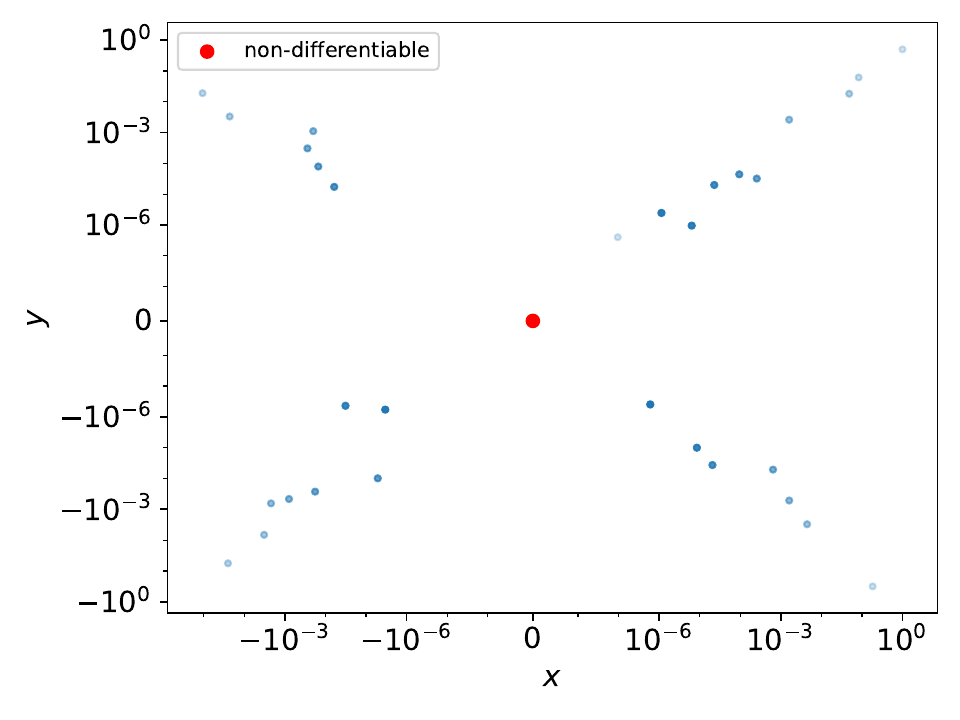}}
\subfloat[]{\includegraphics[width=1.7in]{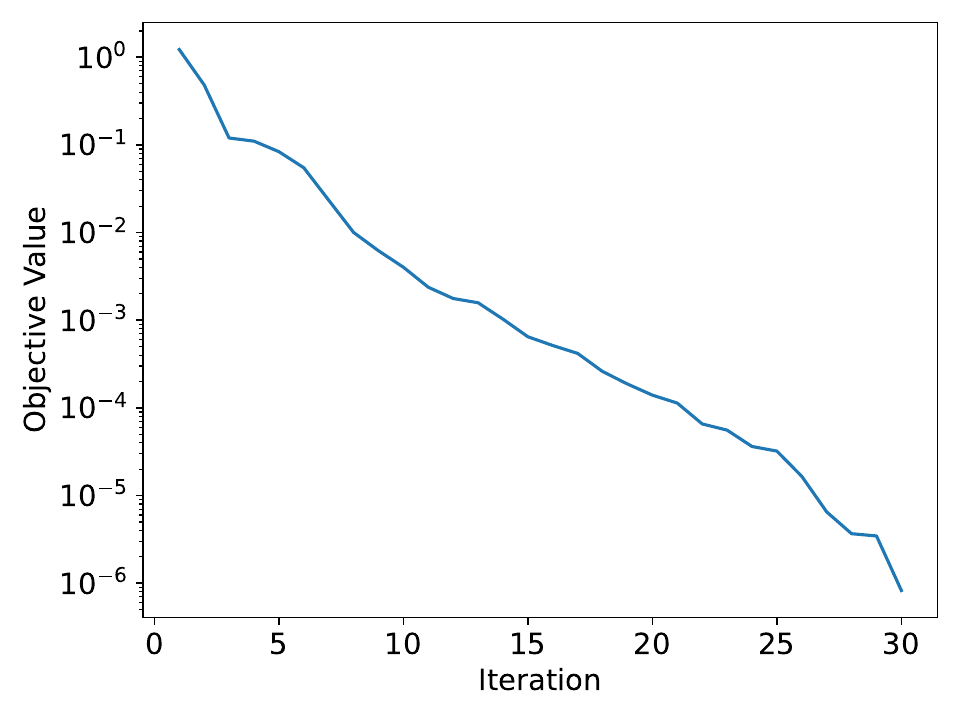}}

\caption{\textbf{\label{fig:norm}Differentiability and Linear Convergence on} \textbf{Nonsmooth}\textbf{ Objectives.} We use BFGS to minimize the nonsmooth objective $h_\mathrm{ellipse}\left(x,y\right)$ that achieves minimum value $0$ at unique minimizer $\bx = \left(0,0\right)$. BFGS is implemented exactly as described in \citet{lewis2009}, initialized at $\left(1,0.5\right)$, and ran for 30 iterations. Panel \textbf{(a)} visualizes $h_\mathrm{ellipse}$. It shows the iterates (blue dots) in the $xy$-plane as well as the origin (red) at which $h_\mathrm{max}$ is non-differentiable. Panels \textbf{(b)} and \textbf{(c)}, respectively, show the gradients at and locations of each iterate with darker colors indicating later iterations. Panel \textbf{(b)} shows gradients forming a convex hull that empirically resembles $\partial h_\mathrm{ellipse}\left(0,0\right)$'s two-dimensional structure. Panel \textbf{(c)} (in symmetric-log scale) shows all iterates lie within differentiable regions of the domain. Panel \textbf{(d)} records the linear convergence of the objective value.} \end{figure}

Similarly, Figure \ref{fig:norm} shows an objective where the subdifferential at its minimizer is a two-dimensional ellipse in the $\R^{2}$ domain\begin{equation}
\partial h_\mathrm{ellipse}\left(\bx\right)=\left\{ \left(x,y\right):\ x^{2}+\frac{1}{2}y^{2}\leq1\right\}.\label{eq:ell_subdiff}
\end{equation}
The ``full-dimensionality'' of $\partial h_\mathrm{ellipse}\left(\bx\right)$ captures the fact that $h_\mathrm{ellipse}$ is non-differentiable in all $\R^{2}$ directions.

\FloatBarrier

Our focus on subdifferentials is also empirically motivated: When employing first-order methods to minimize some objective $h$, the gradients of iterates near the global minimizer $\bx$ typically have a convex hull whose dimension empirically coincides with the dimension of the subdifferential $\partial h\left(\bx\right)$. Figures \ref{fig:phenom}a,b and \ref{fig:norm}a,b demonstrate this behavior. Moreover, this phenomenon occurs despite none of the iterates actually achieving the global minimum (Figures \ref{fig:phenom}c and \ref{fig:norm}c). This is unsurprising since $\partial h (\bx)$ is the convex hull of the limit of gradients at nearby points \citep[Theorem 2.5.1]{clarke1990optimization}. This motivates the creation of a ``survey'' of points near $\bx$ such that their gradients capture the dimensionality of $\partial h (\bx)$. Introducing the max-of-smooth model clarifies this intuition.

\subsubsection{Max-of-Smooth Model}\label{sec:MOSintro}

As is well-known, Fenchel conjugacy allows us to represent any continuous, convex function as the supremum of a family of affine functions (see, for example, \citet[Chapter 3.3]{borwein2010convex}). More generally, continuous, convex functions are instances of \textit{lower-$\C^1$ functions}, which are functions representable as maximums of smooth and compactly parameterized functions.

\begin{defn}
\textbf{\label{def:LC1} (Lower-$\C^{1}$ Functions; \citet{spingarn1981}, \citet[Cor. 3]{daniilidis2004})}

\noindent A locally Lipschitz function $h:\R^{n}\to\R$ is called \textit{lower-}$\C^{1}$ if, for every $\bx\in\R^{n}$, there exists a compact parameterizing set $\T$, a neighborhood $\X$ around $\bx$, and functions $g:\T\times\X\to\R$ such that $g\left(t,x\right)$ and $\nabla_{x}g\left(t,x\right)$ are jointly continuous in $t$ and $x$; and \begin{equation}
h\left(x\right)=\max_{t\in\T}g\left(t,x\right).\label{eq:LC1Max}
\end{equation}
Lower-$\C^{1}$ functions are exactly the locally-Lipschitz and approximately convex functions, where ``approximately convex'' relaxes the canonical notion of convexity and is defined in \citet{daniilidis2004}.
\end{defn}

Lower-$\C^{1}$ functions subsume the entire class of continuous, convex objectives. Definition \ref{def:LC1} implies that, for any global minimizer $\bx$ of the objective $h$, there exists a family of functions $\left\{ g\left(t,\cdot\right)\right\} _{t\in\T}$ satisfying \eqref{eq:LC1Max} at $x{=}\bx$. This family also allows us to describe $h$'s subdifferential at its minimizer:\begin{equation}
\partial h\left(\bx\right)=\Conv\left\{ \nabla_{x}g\left(t,\bx\right):\ t \in \T^\prime\right\},\label{eq:subdiff_max}
\end{equation}
where $\T^\prime = \left\{t:\ g\left(t,\bx\right)=h(\bx)\right\}$. For more details on \eqref{eq:subdiff_max}, see \citet[Theorem 10.31]{rockafellar2009variational}. This characterization inspires us to consider a structurally revealing approximating model for convex, nonsmooth objectives where $\T$ is \textit{finite} and $\T = \T^\prime$. In this setting, \eqref{eq:LC1Max} reduces to a \textit{max-of-smooth} function:\begin{equation}
f\left(x\right)=\max_{i=1,\dots,k}f_{i}\left(x\right),\label{eq:max_mod}
\end{equation}
where $k$ is a finite integer and each $f_{i}$ is a differentiable function. Just like in classical analyses of GD (see, for example, \citet[Chapter 10.6]{beck2017}), we will assume that each $f_i$ is a $\C^2$ function and strongly convex.

When the objective is a max-of-smooth function \eqref{eq:max_mod}, we could ``survey'' the landscape of $f$ by somehow obtaining a point, $s_i$, from each region $\cR_i \equiv \left\{x: \ f_i(x) > f_j(x) \ \forall j\neq i \right\}$ for all $i$. If the gradients $\left\{\nabla f_i(\bx)\right\}_{i=1}^k$ are moreover affinely independent---meaning $\bx$ is a ``nondegenerate'' minimizer---and the $s_i$'s are sufficiently close to $\bx$, the dimension of the convex hull of the survey gradients, $\Conv\left(\left\{\nabla f_i(s_i)\right\}_{i=1}^k\right)$, will  coincide with the dimension of the objective's subdifferential at its minimizer, $\partial f(\bx)=\Conv\left(\left\{\nabla f_i(\bx)\right\}_{i=1}^k\right)$. Moreover, the dimension of $\partial f (\bx)$ is then $k{-}1$. In other words, it is exactly one less than the size of the survey.

\begin{SCfigure}[2][h]\centering
\includegraphics[width=2.1in]{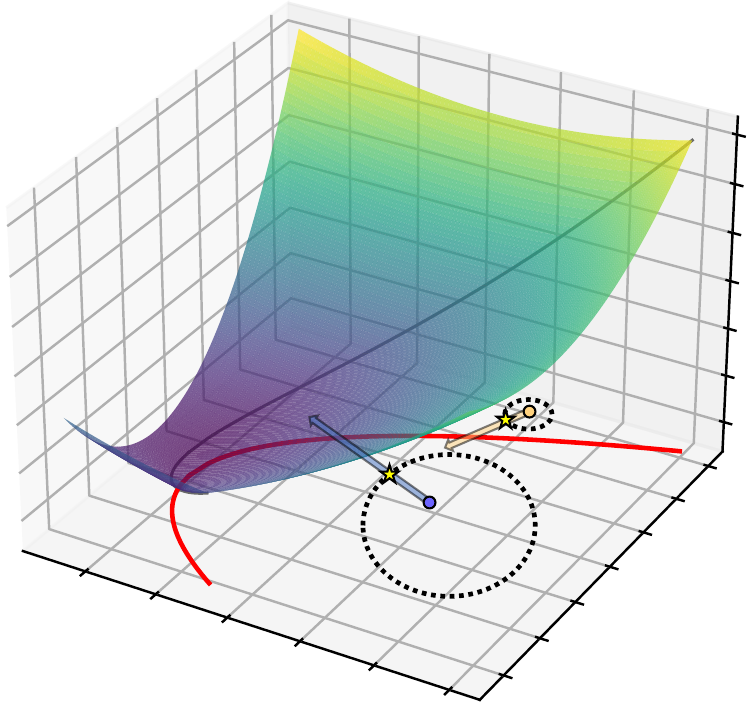}

\caption{\textbf{\label{fig:SDintuition} Intuition of Survey Descent.} Abstract depiction of Survey Descent on max-of-smooth objectives \eqref{eq:max_mod} with $k{=}2$ components and using two survey points (circular dots). Each survey point, $s_i$, is associated with one $f_i$-component of the objective for $i{\in}\{1,2\}$. In the $i$-th Survey Descent subproblem, a gradient step (arrow) from the $i$-th survey point is projected onto the subproblem's constraints (illustrated with  dashed-line boundaries). The constraints \eqref{eq:SD_constraints} prevent subproblem outputs (stars) from crossing the ``nonsmooth boundary'' (red line). In other words, the $i$-th output must remain in the subregion of the domain where the $f_i$-component is ``active''.} 
\end{SCfigure}

Next, observe from Definition \ref{def:SD_Method} that each Survey Descent subproblem simply performs one projected-GD step onto the feasible region defined by \eqref{eq:SD_constraints}. The max-of-smooth objective intuitively motivates this region. In this case, when Survey Descent is equipped with an initializing survey consisting of one point in each $\cR_i$ as discussed above, $f\left(s_{i}\right)=f_{i}\left(s_{i}\right)$ and $\nabla f\left(s_{i}\right)=\nabla f_{i}\left(s_{i}\right)$ for all $i$. Then, since the components $f_{i}$ are $L$-smooth, the constraints for the $i$-th Survey Descent subproblem \eqref{eq:SD_constraints}  restrict solutions to a region where the linear lower bound of $f_{i}$ is at least the quadratic upper bounds of the remaining $f_{j}$, $j\neq i$. Therefore, when the $i$-th subproblem is feasible, its output $s_{i}^{+}$ would necessarily remain within $\cR_{i}$. Section \ref{sec:MOS} formalizes these ideas in detail, but their geometry is simple: Figure \ref{fig:SD_example} shows the behavior of Survey Descent on a simple max-of-smooth objective on $\R^{2}$ while Figure \ref{fig:SDintuition} presents an abstract illustration of the above intuitions.
\FloatBarrier
\FloatBarrier

Despite our discussions of $\left\{f_i\right\}_{i=1}^k$, note the very important feature that \textit{performing the Survey Descent iteration on max-of-smooth function objectives does not require access to individual $f_i$'s}. It only involves function-value and gradient evaluations of the overall objective $f$. Moreover---although Survey Descent is intuitively motivated by the max-of-smooth function objective and possesses a structurally revealing local linear convergence theory on these objectives (Sections \ref{sec:MOS}-\ref{sec:lconv})---it is broadly implementable on \textit{any} objective function. Indeed, Survey Descent still displays promising behavior even on objectives \textit{not} expressible as the maximum of finitely-many smooth functions. For instance, Figures \ref{fig:non_max} and \ref{fig:non_max_big} show the visibly linear convergence of Survey Descent on max-matrix-eigenvalue objectives whose eigenvalues are nonsmooth functions of the input variables.

\subsection{Related Works}\label{sec:related}

A well-known example of \citet{nemirovski1983} shows sublinear worst-case global complexity for any first-order method for arbitrary nonsmooth convex objectives.  That example focuses on an initial sequence of iterates of length less than the dimension of the domain. In contrast, we prove that Survey Descent achieves linear convergence, but only locally and asymptotically, and on an interesting subclass of objectives.

Also worth noting for comparison is classical work on minimizing  the maximum of finitely-many \textit{given} smooth functions $f_i$, where one can access each component $f_i$ rather than just their pointwise maximum, as in Survey Descent. Such problems are standard in classical optimization, easily solved by reduction to classical nonlinear programs.  \citet{kiwiel1983phase} developed a more sophisticated two-phase method, first identifying a search direction by solving a quadratic program built from the function values $f_i$ and gradients $\nabla f_i$, and then performing a line search. More recently, \citet[Scheme 2.3.13]{nesterov2003introductory} computes a linearization of each $f_i$-component, $\left\{\ell^{f_i}_{\tilde{x}}(x)\right\}_{i=1}^k$, at some iterate $\tilde{x}$; performs a proximal-point iteration on $\max_i \ell^{f_i}_{\tilde{x}}(x)$; and then determines the next iteration using a carefully-designed momentum step. Both \citet{kiwiel1983phase} and \citet[Scheme 2.3.13]{nesterov2003introductory} achieve global, linear convergence on strongly convex, max-of-smooth function objectives. 

At large scale, motivated by modern ML applications, \citet{carmon2021} describes a more sophisticated approach. Using function-value and gradient evaluations of the $f_i$-components, \citet{carmon2021} implements a Ball Regularized Optimization Oracle (BROO) that returns the minimizer of the objective---subject to a proximal-regularizer---within a ball around a queried point. Using this BROO, \citet{carmon2021} then designs an efficient minimization procedure by repeatedly updating an iterate using subroutines of line searches and BROO calls. \citet{carmon2021} prove global, sublinear convergence rates for this procedure on max-of-smooth function objectives. 

Within the derivative-free setting, researchers have also studied max-of-smooth objectives \citep{hare2013derivative} as well as generalizations into (potentially nonconvex) objectives that are the composition of multiple component functions \citep{khan2018manifold}. By leveraging careful line-searches and (approximations to) active set information at individual iterates, these methods guarantee convergence using only function-value evaluations of individual objective components without any need for gradient evaluations. Convergence rates are challenging to derive in derivative-free settings and are generally sublinear where they do exist (for example, in \citet{garmanjani2016trust}). For more detailed discussions on works investigating such derivative-free methods, see \citet{khan2018manifold}.

\citet{kiwiel1983phase}, \citet{nesterov2003introductory}, \citet{carmon2021}, and the derivative-free methods above differ from Survey Descent, since they crucially rely on access to the individual components of a max-function objective. In comparison, Survey Descent only assumes a first-order oracle that returns the function value and gradient of $f$ \textit{without} access to any component $f_{i}$'s. Thus, Survey Descent is implementable on \textit{any} objective, even those without max-of-smooth structure, suggesting a promising future research direction. Also notable is the distinction in computational style:  Survey Descent subproblems consist of $k$ parallel projected-GD steps, updating an entire collection of survey points rather than only one iterate as in the works above.

On the other hand, more sophisticated structural oracles can lead even to superlinear convergence on nonsmooth objectives, as remarked in several prior works. For eigenvalue optimization, \citet{oustry2000second} and \citet{noll2005spectral} develop quadratically convergent procedures by incorporating in particular a Hessian oracle within a nearby $\mathcal{U}$-subspace \citep{lemarechal2000U}. On general nonsmooth objectives, using $\mathcal{U}$-Hessians and proximal operators on the objective, \cite{mifflin2005algorithm} present a $\mathcal{VU}$-algorithm whose serious steps converge superlinearly. In contrast, although Survey Descent conceptually relies on a max-of-smooth model, its implementation is general, requiring only function and gradient evaluations.

\sloppy Among multipoint methods, one popular class for nonsmooth optimization are bundle methods \citep{lemarechal1975extension,wolfe1975method}, which possess a long, successful history (described in \citet{mifflin2012science}). In its most transparent proximal form, they rely on a multipoint collection (a ``bundle'') to build a piecewise-linear, cutting-plane model that is minimized in each step with a proximal operator around a current iterate (a ``center''). If the proximal-step outputs a point that sufficiently decreases the function-value, the method takes a ``serious step'' updating the center with the outputted point; otherwise, the method takes a ``null step'' which does not update the center but adds the outputted point to the running point collection to improve the cutting-plane model. 

The related ``level-bundle'' approach  (\citet{lemarechal1995new,kiwiel1995proximal,de2014level}) also shares some similarities with Survey Descent. Each iteration projects the current center onto a sublevel set for the current cutting plane model, the update either being accepted if the objective value decreases satisfactorily---a serious step---or rejected, in which case the cutting plane model is updated. The projection ingredient is similar in both methods, but Survey Descent solves $k$ \textit{parallel} subproblems at each iteration whereas level-bundle methods perform projections sequentially, enforcing objective decrease for serious steps.

Bundle methods work well in practice \citep{borghetti2003lagrangian,emiel2010incremental,sagastizabal2012divide,van2014joint} and enjoy a robust global convergence theory:  see \citet{oliveira2014bundle} for a comprehensive survey. In particular, relative to the number of {\em serious steps}, bundle methods converge linearly on convex objectives whose subdifferentials satisfy certain growth conditions away from the minimizer, as a consequence of a more general framework of \citet{robinson1999linear} and developed further in \citet{atenas2021unified}.\footnote{Two independent related works \citep{atenas2021unified,davis2022nearly} were announced recently.  The first is a flexible, unified analysis of linearly convergent descent methods on weakly convex objectives, in particular covering the serious step sequence for bundle methods.  A transparent analysis of Survey Descent in this framework is not immediately apparent;  the original approach we present here is direct.  The second announcement presents an interesting blackbox randomized first-order method that is nearly linear convergent with high probability.}

However, relative to both null and serious steps, prior published literature \citep{kiwiel2000,du2017,diaz2021} have only derived sublinear convergence guarantees, even when the objective is strongly convex (although \citet{davis2022nearly} is a recent promising advance). Two recent works \citep{liang2021proximal,liang2021unified} unify and present optimal iteration complexities for proximal bundle methods, all of which converge sublinearly, even on strongly convex objectives. In contrast, we show in this paper that the Survey Descent iteration---at least in the case of strongly convex, max-of-smooth function objectives---achieves a \textit{local}, \textit{linear} convergence rate.

Finally, we remark that Survey Descent surveys and bundle method bundles serve different functional roles. In bundle methods, null steps query points sequentially to improve the objective model around the center. In Survey Descent, there are neither auxiliary null steps nor  a center. Instead, a fixed-size survey of points is cautiously updated \textit{in parallel} trying to mimic the steady progress of GD.

\section{Survey Descent on the Max-of-Smooth Objectives}\label{sec:MOS}

We are particularly interested in convex objective functions $f$ with a \textit{nondegenerate minimizer} $\bx$:\begin{equation}
    \bx = \arg\min_{x \in \R^n} f(x) \text{ and } 0\in\Relint\left( \partial f(\bx) \right),\label{eq:nondegen}
\end{equation}
where $\Relint(\cdot)$ denotes the relative interior. We further assume the following structure around $\bx$.
\begin{defn}\label{def:strong_max}
A function $f:\R^n \to \R$ is a \textit{strong $\C^2$ max function} if it is locally expressible near $\bx$ as\begin{equation}\label{eq:f_objective}
f\left(x\right)=\underset{i=1,...k}{\max}f_{i}\left(x\right),
\end{equation}
where $k$ is some finite number, the \textit{components} $\left\{f_i\right\}_{i=1}^k$ are $\C^2$-functions satisfying $f_i(\bx)=f(\bx)$ for all $i$, their gradients $\left\{\nabla f_i(\bx)\right\}_{i=1}^k$ are affinely independent, and their Hessians $\nabla^2 f_i (\bx)$ are positive definite. As a consequence, there exists constants $\delta,L>0$ such that their Hessians satisfy\begin{equation}
    \delta I \preceq \nabla^2 f_i(x) \preceq L I \ \forall \ i=1,\dots,k, \label{eq:hessian}
\end{equation} 
for $x$ near $\bx$---where $I$ is the $n{\times}n$ identity matrix.
\end{defn}

Our analysis in this paper is entirely local, but we \underline{assume for the rest of this paper} that $f$ is a strong $\C^2$ max function with nondegenerate minimizer $\bx$. For simplicity, we also assume that the properties in Definition \ref{def:strong_max} hold \textit{globally} for all $x\in \R^n$. In this setting, we refer to $k$---which is necessarily unique---as the \textit{degree} of $f$. Moreover, note that \eqref{eq:hessian} implies all components $\left\{f_i\right\}_{i=1}^k$ are $L$-smooth, both $f$ and its components are $\delta$-strongly convex, and $\bx$ is $f$'s \textit{unique} minimizer. 

To study Survey Descent (Definition \ref{def:SD_Method}) on the objective $f$, we choose the above-presented $L$ as the step-control parameter and define the following notions of \textit{valid} and \textit{minimizing} surveys.

\begin{defn}
\textbf{\label{def:validS} ($\S, \bX$-Surveys and Validity)} For the objective $f$, define a \textit{survey} as a matrix $\S \equiv \left[s_{1}{,}\dots{,}s_{k}\right] \in \mathbb{\R}^{n\times k}$ consisting of $k$ columns\footnote{We can equivalently consider $\S$ as a set $\{s_i\}_{i=1}^k$. We will use the matrix and set interpretations interchangeably.} (the ``survey points''), where $k$ is the degree of $f$; \textit{valid surveys} as surveys satisfying $f_{i}\left(s_{i}\right){>}f_{j}\left(s_{i}\right)$ for all $i{\neq}j$; and the \textit{minimizing survey} as $\bX \equiv [\bx{,}\dots{,}\bx]$ whose survey points are all equal to $\bx$.
\end{defn}
On the space of surveys, we adopt the norm\begin{equation}
\left\| \S\right\| \equiv\left\| \left[s_{1},\dots,s_{k}\right]\right\| _{2,\infty}=\max_{i=1,\dots,k}\left\| s_{i}\right\| _{2},\label{eq:surv_norm}
\end{equation}
which allows us to quantify the distance between a survey $\S$ from $\bX$ using $\left\| \S{-}\bX \right\|$, an intuitive measure of the distance to the furthest survey point.  The precise choice of norm is immaterial to our computations and theory.

Given a valid $\S$ sufficiently close to $\bX$, Survey Descent exhibits many desirable properties when minimizing $f$. As a preliminary, note that \eqref{eq:nondegen} and Definition \ref{def:strong_max} implies there exists unique \textit{critical weights} $\left\{ \blam_{i}\right\}_{i=1}^{k}{>}0$ such that 
\begin{align}
\sum_{i=1}^{k}\blam_{i} & =1,\label{eq:lam_sum1}\\
\sum_{i=1}^{k}\blam_{i}\nabla f_{i}(\bx) & =0.\label{eq:lam_sum0}
\end{align}
As we will see, the Lagrange multipliers of Survey Descent subproblems will lie ``close'' to $\left\{ \blam_{i}\right\}_{i=1}^{k}{>}0$, while the subproblem outputs $\left\{s_i^+\right\}_{i=1}^k$ will lie near  $\bx$. To characterize the order-of-magnitude of these distances, we adopt the ``Big-Oh'' notation. In particular, for a mapping $g:\mathbf{E}{\to}\F$ between two Euclidean spaces and letting $\left\| \cdot\right\|^p$ denote an arbitrary Euclidean norm raised to the $p$-th power, we use the notation 
$g\left(x\right)=O\left(\left\| x\right\| ^{p}\right)$
to indicate the property that there exists a constant $K{>}0$ such that $ \left\| g\left(x\right)\right\|{\leq}K\left\| x\right\| ^{p}$ holds for all small $x$. By itself, we let $\mathrm{O}\left(\left\| x\right\| ^{p}\right)$ denote an element of the class of all functions with this property.

With this terminology, we present our primary theoretical setting as well as our first result.

\begin{setting} \label{setting:localSD} \textbf{(Local Analysis of Survey Descent)}
Consider an iteration of Survey Descent on a strong $\C^2$ max function objective, $f$, with a survey 
$\S$ that is valid and sufficiently close to $\bX$. For $i{=}1{,}\dots{,}k$, denote the output of the $i$-th Survey Descent subproblem by $s_i^+$.
\end{setting}

\begin{restatable}{thm}{PBmbp}\label{thm:SD_MBP} \textbf{(Local Feasibility, Uniqueness, Tightness, and Smoothness of Survey Descent)} Assume Setting \ref{setting:localSD}. Then for all $i$, the Survey Descent subproblem $\left(P_{i}^{\S}\right)$ is feasible and has a unique solution $s_{i}^{+}$, which satisfies all the constraints with equality, with unique associated Lagrange multipliers $\lambda_{j}^i$ \edit{for the $j$-th constraint within the $i$-th subproblem}  (for $j\neq i$). The solutions and multipliers depend smoothly on the input survey $\S$, and satisfy\begin{equation}
s_{i}^{+}\left(\S\right) ~=~\bx+\mathrm{O}\left(\left\| \S-\bX\right\| \right) \qquad \text{and} \qquad \lambda_{j}^i\left(\S\right) ~=~ \blam_{j}+\mathrm{O}\left(\left\| \S-\bX\right\| \right)\ \forall\ j\neq i.\label{eq:SD_forms}
\end{equation}
\end{restatable}
\begin{proof}
Follows from a routine analysis of the subproblem's first-order optimality conditions combined with the implicit function theorem.
\end{proof}
We refer to Appendix B \ifrevise of \citet{han2021survey} \fi for a more direct and elementary proof that is computationally revealing. In particular, it shows we can compute $s_{i}^{+}\left(\S\right)$ and $\lambda_{j}^i\left(\S\right)$ entirely with linear algebra routines and a scalar square-root.

Theorem \ref{thm:SD_MBP} plays a key role in deducing many theoretical results in this paper. Most immediately, under the assumptions of the theorem, the output survey $\S^+=\left\{s_i^+\right\}_{i=1}^k$ of a Survey Descent iteration will also be valid. First, note that $f$ is differentiable at every valid survey point and the associated gradients satisfy\begin{equation}\label{eq:valS}
\nabla f\left(s_{i}\right)=\nabla f_{i}\left(s_{i}\right) \ \forall \ i=1,\dots,k,
\end{equation}
which leads to the following observation.
\begin{observe}
\textbf{(Survey Descent Iterations on Max-of-Smooth Objectives)\label{obs:SD_MOS}} Assume that $\S$ is a valid survey for the objective $f$. Then, for all $i$, the $i$-th subproblem $\left(P_{i}^{\S}\right)$ of Survey Descent (Definition \ref{def:SD_Method}) is equivalent to \begin{equation}\begin{aligned}\min_{x}\ & \left\| x-\left(s_{i}-\frac{1}{L}\nabla f_{i}\left(s_{i}\right)\right)\right\| _{2}^{2}\\
\mathrm{s.t.}\ & \ell_{s_j}^{f_j}(x)+\frac{L}{2}\left\| x-s_{j}\right\| _{2}^{2}\leq \ell_{s_i}^{f_i}(x)\ \forall\ j\neq i.
\end{aligned}
\label{eq:raw_form_analytic}
\end{equation}
\end{observe}

Using Theorem \ref{thm:SD_MBP} and Observation \ref{obs:SD_MOS}, we deduce the aforementioned validity of $\S^+$.

\begin{thm} \label{thm:strictval} \textbf{(Preservation of Validity)} Assuming Setting \ref{setting:localSD}, the output survey $\S^{+}$ of a Survey Descent iteration is valid.
\end{thm}
\begin{proof}
Consider any fixed $i$. If $s_{i}^{+}=s_{i}$, then the fact
that $f_{i}\left(s_{i}^{+}\right){>}f_{j}\left(s_{i}^{+}\right)$ for all $j\neq i$ immediately follows since Theorem \ref{thm:SD_MBP} assumes $\S$ is valid.

Now, consider the $s_{i}^{+}{\neq}s_{i}$ case. By Theorem \ref{thm:SD_MBP}, $s_{i}^{+}$ satisfies the constraints of $\left(P_{i}^{\S}\right)$ with equality. Using Observation \ref{obs:SD_MOS}, we express this as\begin{equation}
\ell_{s_i}^{f_i}(s_i^+) = \ell_{s_j}^{f_j}(s_i^+)+\frac{L}{2}\left\| s_{i}^{+}-s_{j}\right\|_{2}^{2}\ \forall\ j\neq i.\label{eq:left_valid}
\end{equation}
The $\delta$-strong convexity of $f_{i}$ implies
\begin{equation*}
f_{i}\left(s_{i}^{+}\right)\geq\ell_{s_i}^{f_i}(s_i^+)+\frac{\delta}{2}\left\| s_{i}^{+}-s_{i}\right\| _{2}^{2}>\ell_{s_i}^{f_i}(s_i^+),
\end{equation*}
where the strict inequality follows from $s_{i}^{+}{\neq}s_{i}$. Next, by the $L$-smoothness of $\left\{f_j\right\}_{j=1}^k$, the right-hand side of \eqref{eq:left_valid} is a $s_j$-centered quadratic upper bound of $f_{j}$ evaluated at $s_i^+$. Thus,
\begin{equation*}
\ell_{s_j}^{f_j}(s_i^+)+\frac{L}{2}\left\| s_{i}^{+}-s_{j}\right\|_{2}^{2} \geq f_{j}\left(s_{i}^{+}\right) \ \forall\ j\neq i.
\end{equation*}

\noindent Therefore, 
\begin{equation*}
f_{i}\left(s_{i}^{+}\right) > \ell_{s_i}^{f_i}(s_i^+) = \ell_{s_j}^{f_j}(s_i^+)+\frac{L}{2}\left\| s_{i}^{+}-s_{j}\right\|_{2}^{2} \geq f_{j}\left(s_{i}^{+}\right) \ \forall\ j\neq i.
\end{equation*}
This completes the proof.
\end{proof}

\section{Connecting Survey Descent and Gradient Descent}
\sloppy We will build a local convergence theory by connecting Survey Descent to projected GD-steps on the smooth components $\left\{f_i\right\}_{i=1}^k$ of the objective $f$. To form this connection, we first identify the affine \textit{$\U$-subspace} possessing the following equivalent characterizations:
\begin{align}
\U & \equiv \bx + \Span\left\{ \nabla f_{i}(\bx)-\nabla f_{j}(\bx):1\leq i,j\leq k\right\} ^{\perp} \label{eq:Udef1}\\
 & = \bx +  \left\{ \sum_{i}^k\gamma_{i}\nabla f_{i}(\bx):\sum_{i}^k\gamma_{i}=0\right\} ^{\perp} \label{eq:Udef2}\\
 & = \bx +  \left\{x\in \R^n: \ \nabla f_i(\bx)^T x = \nabla f_j(\bx)^T x \ \forall \ 1\leq i,j\leq k\right\}.\label{eq:Udef3}
\end{align}
$\U$ is so-named because the $\bx$-centered, linear subspace $\left\{\U -\bx\right\}$ is commonly called the ``$\mathcal{U}$-subspace'' \citep{lemarechal2000U} and captures the directions in which $f$ is smooth around $\bx$. Theorem \ref{thm:approx_s} will show that the outputs of Survey Descent subproblems are approximately---up to an $\mathrm{O}\left(\|\S - \bX\|^2\right)$ term---a convex combination of $\U$-projected GD-steps on the smooth components $\left\{f_i\right\}_{i=1}^k$. 

To prove Lemmas \ref{lem:r_grad}-\ref{lem:Up_component} leading to Theorem \ref{thm:approx_s}, we adopt the following notation for the gradients and Hessians of the components of $f$ (Definition \ref{def:strong_max}):\begin{equation}
    a_i \equiv \nabla f_i(\bx) \text{ and } A_i \equiv \nabla^2 f_i(\bx) \ \forall \  i=1,\dots,k. \label{eq:A_notation}
\end{equation}
Additionally, assuming $\bx{=}0$ without loss of generality will simplify many of our proofs. In this case, the second-order Taylor expansion of $f_i$ around zero is\begin{equation}
    f_{i}\left(x\right)= f(0) + a_{i}^{T}x+\frac{1}{2}x^{T}A_{i}x+r_{i}\left(x\right) \ \forall \ i=1,\dots,k,\label{eq:taylor}
\end{equation}
where $r_i(x)$ is a residual function.

\begin{lem}
\label{lem:r_grad} For all $i$, the residual function $r_i(x)$ in \eqref{eq:taylor} satisfies $r_i(x)=\mathrm{O}\left(\|x\|^2\right)$ and $\nabla r_{i}\left(x\right)=\mathrm{O}\left(\left\| x\right\|\right)$.
\end{lem}
\begin{proof}
$r_i(x)=\mathrm{O}\left(\|x\|^2\right)$ follows directly from Taylor's Theorem. Next, observe that 
\begin{equation*}
    r_{i}\left(x\right)=f_{i}\left(x\right)-\left(f(0) + a_{i}^{T}x+\frac{1}{2}x^{T}A_{i}x\right).
\end{equation*}
Thus, since $f_i$ is $\C^2$, the residual $r_i$ must also be $\C^2$, so its gradient and Hessian are well-defined. This implies $\nabla r_{i}\left(x\right)$ is then $\C^{1}$. Taylor expanding $\nabla r_i$ around zero then gives the desired $\nabla r_{i}\left(x\right)= \mathrm{O}\left(\| x\|\right)$
after observing that $\nabla r_i(0)$ is the zero-vector and $\nabla^2 r_i(0)$ is the zero-matrix (since $r_i$ is the residual of the Taylor expansion of $f_i$ around zero).
\end{proof}

We now prove and present the main results of this section.

\begin{lem}\label{lem:U_component} Assume Setting \ref{setting:localSD} and suppose $\bx{=}0$. Then, $s_{i}^{+}$ satisfies
\begin{equation*}
\proj\left[s_{i}^{+}\right]=\frac{1}{2-\blam_{i}}\proj\left[s_{i}-\frac{\blam_{i}}{L}\nabla f_{i}\left(s_{i}\right)\right]+\sum_{j\neq i}\frac{\blam_{j}}{2-\blam_{i}}\proj\left[s_j - \frac{1}{L} \nabla f_j(s_j)\right]+\mathrm{O}\left(\left\| \S\right\| ^{2}\right) \ \forall \ i=1,\dots,k.
\end{equation*}
where $\proj$ is the projection operator onto the $\U$-subspace \eqref{eq:Udef1}.
\end{lem}
\begin{proof}
Without loss of generality, assume $i=k$. Scale the objective of problem \eqref{eq:raw_form_analytic} by $\frac{L}{2}$ and consider the equivalent problem
\begin{equation}\begin{aligned}\min_{x}\ & \frac{L}{2} \left\| x-\left(s_{k}-\frac{1}{L}\nabla f_{k}\left(s_{k}\right)\right)\right\| _{2}^{2}\\
\mathrm{s.t.}\ & \ell_{s_j}^{f_j}(x)+\frac{L}{2}\left\| x-s_{j}\right\| _{2}^{2}\leq \ell_{s_k}^{f_k}(x)\ \forall\ j\neq k.
\end{aligned}\label{eq:scaledPk}
\end{equation}
By Theorem \ref{thm:SD_MBP} and Observation \ref{obs:SD_MOS}, there exists unique Lagrange multipliers $\left\{\lambda_j\right\}_{j=1}^{k-1}$ (we omit the $k$-superscript for simplicity) satisfying the following first-order optimality condition for \eqref{eq:scaledPk}:
$$
    0=L\left(s_{k}^{+}-\left(s_{k}-\frac{1}{L} \nabla f_{k}\left(s_{k}\right)\right)\right)+\sum_{j=1}^{k-1} \lambda_{j}\left(\edit{\nabla f_{j}\left(s_{j}\right)}+L\left(s_{k}^{+}-s_{j}\right)-\nabla f_{k}\left(s_{k}\right)\right)
    $$
After rearranging this gives\begin{equation}
0= L\left(1+\sum_{j=1}^{k-1}\lambda_{j}\right)s_{k}^{+}-L\left(s_{k} + \sum_{j=1}^{k-1} \lambda_{j} s_{j}\right)+\nabla f_k(s_k) +\sum_{j=1}^{k-1}\lambda_{j}\left(\nabla f_j(s_j) - \nabla f_k(s_k)\right).
\label{eq:FOC2}
\end{equation}
Define the following scalars:
\begin{align}
    \lambda_{k} & \equiv  1-\sum_{j=1}^{k-1}\lambda_{j},\label{eq:lambda}\\
    \mu_{j} & \equiv  \blam_{j}-\lambda_{j} \ \forall \ j=1,\dots,k. \label{eq:mu}
\end{align}
For $j\neq k$, Theorem \ref{thm:SD_MBP} gives
$\mu_{j}=\mathrm{O}\left(\left\| \S\right\| \right)$ .
For $j=k$, combining \eqref{eq:lambda}-\eqref{eq:mu} with equation \eqref{eq:lam_sum1} implies
\begin{equation*}
\mu_{k} =\blam_{k}-\left(1-\sum_{j=1}^{k-1}\lambda_{j}\right) =\blam_{k}-\left(1-\sum_{j=1}^{k-1}\blam_{j}\right)+\mathrm{O}\left(\left\| \S\right\| \right) =\mathrm{O}\left(\left\| \S\right\| \right).
\end{equation*}
Thus,\begin{equation}
    \mu_j = \mathrm{O}(\|\S\|) \ \forall \ j=1,\dots,k.\label{eq:mu_size}
\end{equation}
Substituting \eqref{eq:mu} into \eqref{eq:FOC2}, applying identities \eqref{eq:lam_sum1}-\eqref{eq:lam_sum0}, collecting $\mathrm{O}\left(\|\S\|^2\right)$ terms, and solving for $s_k^{+}$ gives\begin{equation}s_{k}^{+}= \frac{1}{2-\blam_{k}}\left[s_{k}-\frac{\blam_{k}}{L}\nabla f_{k}\left(s_{k}\right)\right]+\sum_{j=1}^{k-1}\frac{\blam_{j}}{2-\blam_{k}}\left[s_j - \frac{1}{L} \nabla f_j(s_j)\right]+\frac{1}{L\left(2-\blam_{k}\right)}\sum_{j=1}^{k}\mu_{j}\nabla f_j(s_j) +\mathrm{O}\left(\left\| \S\right\| ^{2}\right).\label{eq:pre_sform}
\end{equation}
Adopting the notation \eqref{eq:A_notation} and applying the Taylor expansion \eqref{eq:taylor}, observe\begin{equation}
    \sum_{j=1}^{k}\mu_{j}\nabla f_j(s_j) = \sum_{j=1}^{k}\mu_{j}\left(a_j + A_j s_j + \nabla r_j (s_j)\right) = \sum_{j=1}^{k}\mu_{j}a_j +\mathrm{O}\left(\left\| \S\right\|^{2}\right),\label{eq:sum_muA}
\end{equation}
where last equality follows from \eqref{eq:mu_size} and Lemma \ref{lem:r_grad}. Substituting \eqref{eq:sum_muA} into \eqref{eq:pre_sform} then leads to\begin{equation}
s_{k}^{+}=\frac{1}{2-\blam_{k}}\left[s_{k}-\frac{\blam_{k}}{L}\nabla f_{k}\left(s_{k}\right)\right]+\sum_{j=1}^{k-1}\frac{\blam_{j}}{2-\blam_{k}}\left[s_j - \frac{1}{L} \nabla f_j(s_j)\right]+\frac{1}{L\left(2-\blam_{k}\right)}\underbrace{\sum_{j=1}^{k}\mu_{j}a_j}_{(\star)} +\mathrm{O}\left(\left\| \S\right\| ^{2}\right).
\label{eq:s_form}
\end{equation}
Next, again combine \eqref{eq:lambda}-\eqref{eq:mu} with equation \eqref{eq:lam_sum1} to deduce
\begin{equation*}
\sum_{j=1}^{k}\mu_{j}=\sum_{j=1}^{k}\left(\blam_{j}-\lambda_{j}\right)=0.
\end{equation*}
Then, it follows from \eqref{eq:Udef2} that $\proj\left[(\star)\right]=0$.
The desired result then follows from applying $\proj$ to \eqref{eq:s_form}.
\end{proof}

\begin{lem}\label{lem:Up_component}
Assume Setting \ref{setting:localSD}, and let $\bx{=}0$. Then, $s_{i}^{+}$ satisfies
\begin{equation*}
\projp\left[s_{i}^{+}\right]=\mathrm{O}\left(\left\| \S\right\| ^{2}\right) \ \forall \ i=1,\dots,k.
\end{equation*}
where $\projp$ is the projection operator onto the orthogonal complement of the $\U$-subspace \eqref{eq:Udef1}.
\end{lem}
\begin{proof}
Consider any fixed $i$. By Theorem \ref{thm:SD_MBP}, the constraints of $\left(P_{i}^{\S}\right)$
are tight. Thus, by Observation \ref{obs:SD_MOS}, $s_{i}^{+}$ satisfies
\begin{equation*}
    f_{j}\left(s_{j}\right)+\nabla f_{j}\left(s_{j}\right)^{T}\left(s_i^+-s_{j}\right)+\frac{L}{2}\left\| s_i^+-s_{j}\right\| _{2}^{2} = f_{i}\left(s_{i}\right)+\nabla f_{i}\left(s_{i}\right)^{T}\left(s_i^+-s_{i}\right) \ \forall \ j\neq i.
\end{equation*}
Using equation \eqref{eq:taylor} in the above gives, for all $j\neq i$,
\begin{multline*}
    f(0) + a_{j}^{T}s_{j}+\frac{1}{2}s_{j}^{T}A_{j}s_{j}+r_{j}\left(s_j\right)+\left(a_{j}+A_{j}s_{j} + \nabla r_j(s_j)\right)^T\left(s_{i}^{+}-s_{j}\right)+\frac{L}{2}\left\Vert s_{i}^{+}-s_{j}\right\Vert _{2}^{2}\\
    =f(0) + a_{i}^{T}s_{i}+\frac{1}{2}s_{i}^{T}A_{i}s_{i}+r_{i}\left(s_i\right)+\left(a_{i}+A_{i}s_{i}+\nabla r_i(s_i)\right)^T\left(s_{i}^{+}-s_{i}\right).
\end{multline*}
After simplification and rearrangement, we obtain\begin{equation}
\left(a_{j}-a_{i}\right)^{T}s_{i}^{+}=-\left(A_{j}s_{j}-A_{i}s_{i}+\nabla r_{j}\left(s_{j}\right)-\nabla r_{i}\left(s_{i}\right)\right)^{T}s_{i}^{+}-\frac{L}{2}\left\| s_{i}^{+}-s_{j}\right\| _{2}^{2}-c_{j}\left(s_{j}\right)+c_{i}\left(s_{i}\right)  \ \forall \ j\neq i, \label{eq:s_constr}
\end{equation}
where we define the functions
\begin{equation*}
c_{j}\left(x\right)\equiv-\frac{1}{2}x^{T}A_{j}x+r_{j}\left(x\right)-\nabla r_{j}\left(x\right)^{T}x \ \forall \ j=1,\dots,i.
\end{equation*}
By Lemma \ref{lem:r_grad}, $c_j(x)$ is $O\left(\|x\|^2\right)$. Since we assume $\bx{=}0$, by Theorem \ref{thm:SD_MBP}, $s_{i}^{+}=\mathrm{O}\left(\left\| \S\right\| \right)$.
We then deduce that all terms on the right-hand side of \eqref{eq:s_constr} are $\mathrm{O}\left(\left\| \S\right\|^2 \right)$:
\begin{equation}
\left(a_{j}-a_{i}\right)^{T}s_{i}^{+}=\mathrm{O}\left(\left\| \S\right\|^2 \right) \ \forall \ j\neq i.\label{eq:projsp}
\end{equation}
By \eqref{eq:Udef1},
\begin{align*}
\U^{\perp} & =\Span\left\{ a_{i}-a_{j}:1\leq i,j\leq k\right\} =\Span\left\{ a_{j}-a_{i}:\ j \neq i \right\}.
\end{align*}
Combined with the affine independence of $\left\{ a_{j}\right\} _{j=1}^{k}$
(Definition \ref{def:strong_max}), the above implies $\left\{ a_{j}-a_{i}\right\} _{j\neq i}$ forms a basis for $\U^{\perp}$. Therefore, we can express the projection onto $\U^{\perp}$ in terms $\left\{ a_{j}-a_{i}\right\} _{j\neq i}$. The result now follows from equation \eqref{eq:projsp} and standard linear algebra.
\end{proof}

\begin{thm}
\label{thm:approx_s} \textbf{(Gradient Descent Approximation of Survey Updates)} Assume Setting \ref{setting:localSD}.  Then the solution of the $i$-th Survey Descent subproblem satisfies, for all $i$, \begin{equation}
s_{i}^{+}= \ts_{i}+\mathrm{O}\left(\left\| \S - \bX\right\| ^{2}\right) \label{eq:SD_iter}
\end{equation}
where we define\begin{equation}
\ts_{i}\equiv\frac{1}{2-\blam_{i}}\proj\left[s_{i}-\frac{\blam_{i}}{L}\nabla f_{i}\left(s_{i}\right)\right]+\sum_{j\neq i}\frac{\blam_{j}}{2-\blam_{i}}\proj\left[s_j - \frac{1}{L} \nabla f_j(s_j)\right], \label{eq:s_tilde}
\end{equation}
and $\proj$ is the projection operator onto the $\U$-subspace \eqref{eq:Udef1}.
\end{thm}
\begin{proof}
Without loss of generality, assume that $\bx{=}0$. Then, the desired result follows from decomposing $s_{k}^{+} = \proj\left[s_{k}^{+}\right] + \projp\left[s_{k}^{+}\right]$ and applying Lemmas \ref{lem:U_component} and \ref{lem:Up_component}. We can deduce the general $\bx{\neq}0$ case with the change of variables $x{\leftarrow}(x-\bx)$.
\end{proof}

\section{Local Linear Convergence of Surveys\label{sec:lconv}}

Observe that the constituent terms of $\ts_i$ in \eqref{eq:s_tilde} are $\U$-projected gradient steps on the $L$-smooth and $\delta$-strongly convex component functions $\left\{f_i\right\}_{i=1}^k$ (Definition \ref{def:strong_max}). Thus, we can apply classical theory on projected-GD described in standard texts such as \citet[Theorem 10.29]{beck2017}. The idea of tracking iterates relative to their projections onto some nearby ``active'' manifold is familiar in nonsmooth optimization. Recent examples include \citet{burke2020gradient}, \citet{duchi2021asymptotic}, and \citet{davis2021subgradient}.

\begin{thm}
\textbf{\label{thm:PR}(Projected-GD Behavior)} Assume that $\S$ is a valid survey for the strong $\C^2$ max function objective $f$. Then, the following holds for all $i$ and any constant $\tL \geq L$:
\begin{align*}
\left\| \proj\left[s_i - \frac{1}{\tL} \nabla f_i(s_i)\right] - \bx\right\| _{2}^{2} &\leq \left(1-\frac{\delta}{\tL}\right)\left\| s_{i} - \bx \right\| _{2}^{2},\\
f_{i}\left(\proj\left[s_i - \frac{1}{\tL} \nabla f_i(s_i)\right]\right) - f_i(\bx) &\leq \frac{\tL}{2}\left(1-\frac{\delta}{\tL}\right)\left\| s_{i} - \bx\right\| ^{2}.
\end{align*}
\end{thm}
\begin{proof}
Consider any fixed $i$. It is easy to check that $\bx = \arg\min_{x \in \U} f_i(x)$ using first-order optimality conditions combined with \eqref{eq:lam_sum1}-\eqref{eq:lam_sum0} and \eqref{eq:Udef3}. The desired result then follows from classical projected-GD convergence results (given in standard texts such as \citet[Theorem 10.29]{beck2017}) after observing that $\U$ is convex and $f_i$ is $L$-smooth and $\delta$-strongly convex.
\end{proof}

To prove local linear convergence, we show below that the survey points contract towards $\bx$. The corresponding contraction ratios will depend on the minimum critical weight defined as follows:\begin{equation}
    \blam_{\min} \equiv \min_{i=1,\dots,k} \blam_{i} > 0. \label{eq:min_lam}
\end{equation}

\begin{lem}
\label{lem:st_bound} Assume Setting \ref{setting:localSD}, and suppose $\bx{=}0$. Then, the following must hold:
\begin{equation*}
\left\| \ts_{i}\right\| _{2}^{2}\leq\left(1-\frac{\blam_{\min}\delta}{L}\right)\max_{j=1,....,k}\left\| s_{j}\right\| _{2}^{2} \ \forall \ i=1,\dots,k,
\end{equation*}
where $\ts_i$ and $\blam_{\min}$ are as defined in equations \eqref{eq:s_tilde} and \eqref{eq:min_lam}, respectively.
\end{lem}
\begin{proof} Consider any fixed $i$. Then,
\begin{align*}
\left\| \ts_{i}\right\| _{2}^{2} & =\left\| \frac{1}{2-\blam_{i}}\proj\left[s_{i}-\frac{\blam_{i}}{L}\nabla f_{i}\left(s_{i}\right)\right]+\sum_{j\neq i}\frac{\blam_{j}}{2-\blam_{i}}\proj\left[s_j - \frac{1}{L} \nabla f_j(s_j)\right] \right\| _{2}^{2}  \text{ (Substitute \eqref{eq:s_tilde})}\\
 & \leq\frac{1}{2-\blam_{i}}\left\| \proj\left[s_{i}-\frac{\blam_{i}}{L}\nabla f_{i}\left(s_{i}\right)\right]\right\| _{2}^{2}+\sum_{j\neq i}\frac{\blam_{j}}{2-\blam_{i}}\left\| \proj\left[s_j - \frac{1}{L} \nabla f_j(s_j)\right]\right\| _{2}^{2}\text{ (Convexity)}\\
 & \leq\frac{1}{2-\blam_{i}}\left(1-\frac{\blam_{i}\delta}{L}\right)\left\| s_{i}\right\| _{2}^{2}+\left(1-\frac{\delta}{L}\right)\sum_{j\neq i}\frac{\blam_{j}}{2-\blam_{i}}\left\| s_{j}\right\| _{2}^{2} \text{ (Theorem \ref{thm:PR} where $\tL=\frac{L}{\blam_i}$)}\\
 & \leq\left(1-\frac{\blam_{i}\delta}{L}\right)\left[\frac{1}{2-\blam_{i}}\left\| s_{i}\right\| _{2}^{2}+\sum_{j\neq i}\frac{\blam_{j}}{2-\blam_{i}}\left\| s_{j}\right\| _{2}^{2}\right]\\
 & \leq\left(1-\frac{\blam_{\min}\delta}{L}\right)\max_{j=1,....,k}\left\| s_{j}\right\| _{2}^{2}.
\end{align*}
This completes the proof.
\end{proof}

\begin{thm} \textbf{(Survey-Norm Contraction)} \label{thm:s_monotone} Assume Setting \ref{setting:localSD}.  Then the following must hold:
\begin{equation*}
\max_{j=1,...,k}\left\| s_{j}^{+} - \bx \right\| _{2}^{2}\leq\left(1-\frac{\blam_{\min}\delta}{L}\right)\max_{j=1,...,k}\left\| s_{j} - \bx \right\| _{2}^{2}+\mathrm{O}\left(\left\| \S - \bX\right\|^{3}\right),
\end{equation*}
where $\blam_{\min}$ is defined as in \eqref{eq:min_lam}.
\end{thm}
\begin{proof}
Without loss of generality, assume $\bx = 0$. Consider any fixed $i$. By Lemma \ref{lem:st_bound}, $\ts_{i}=\mathrm{O}\left(\left\| \S\right\| \right)$. Then,\begin{align*}
\left\| s_{i}^{+}\right\| _{2}^{2} & =\left\| \ts_{i}+\mathrm{O}\left(\left\| \S\right\| ^{2}\right)\right\| _{2}^{2} \text{ (Theorem \ref{thm:approx_s})}\\
 & =\left\| \ts_{i}\right\| _{2}^{2}+2\left\langle \ts_{i},\mathrm{O}\left(\left\| \S\right\| ^{2}\right)\right\rangle +\mathrm{O}\left(\left\| \S\right\| ^{4}\right)\\
 & \leq\left(1-\frac{\blam_{\min}\delta}{L}\right)\max_{j=1,....,k}\left\| s_{j}\right\| _{2}^{2}+\mathrm{O}\left(\left\| \S\right\| ^{3}\right) \text{ (Lemma \ref{lem:st_bound})}.
\end{align*}
The right-hand side is independent of $i$, so taking the max over $i$ on the left-hand side leads to the desired result. We can deduce the general $\bx{\neq}0$ case with the change of variables $x{\leftarrow}(x-\bx)$.
\end{proof}

Theorem \ref{thm:s_monotone} implies that for input surveys $\S$, sufficiently close to $\bX$, Survey Descent outputs a survey $\S^{+}$ that is strictly \textit{closer} to $\bX$. As a result, the feasibility and validity guarantees of Theorems \ref{thm:SD_MBP} and \ref{thm:strictval} will continue to hold in repeated applications of Survey Descent iterations allowing us to repeat the procedure indefinitely. Letting $\nonneg$ denote the non-negative natural numbers, we formalize these ideas below.

\begin{procedure}\textbf{(Repeated Iterations of Survey Descent) }\label{procedure:SD}
For the strong $\C^2$ max function objective $f$, assume $\S^{0}$ is a valid initial survey sufficiently close to $\bX$, and initialize $t=0$. Iterate the following steps:
\begin{enumerate}
\item Solve the Survey Descent iteration $\left\{ \left(P_{i}^{\S^{t}}\right)\right\} _{i=1}^{k}$
and denote the output $\S^{t+1}\equiv\left(\S^{t}\right)^{+}$.
\item Increment $t\leftarrow t+1$. 
\end{enumerate}
\end{procedure}

\begin{cor}\textbf{(Well-Definedness of Survey Descent Repetitions)}\label{cor:viable} Procedure \ref{procedure:SD} is well-defined:  for each $t=0,1,2\ldots$, the Survey Descent iteration $\left\{ \left(P_{i}^{\S^{t}}\right)\right\} _{i=1}^{k}$ is feasible and produces a valid output survey.
\end{cor}
\begin{proof}
For $\S^{0}$ sufficiently close to $\bX$, Theorem \ref{thm:s_monotone} inductively implies that $\left\| \S^{t+1} - \bX \right\| <\left\| \S^{t} - \bX\right\| $ for all $t\in\nonneg$. In other words, output surveys will always move closer to $\bX$, and Theorem \ref{thm:SD_MBP} and \ref{thm:strictval} continue to hold for Survey Descent iterations $\left\{\left(P_{i}^{\S^{t}}\right)\right\}_{i=1}^{k}$ for all $t\in \nonneg$. Thus, Theorem \ref{thm:SD_MBP} implies the feasibility all Survey Descent iterations while Theorem \ref{thm:strictval} implies the validity of output surveys.
\end{proof}

\begin{observe}\label{obs:allT}
The proof of Corollary \ref{cor:viable} more generally shows that, for all $t \in \nonneg$, all surveys $\S^{t}$ will remain sufficiently close to $\bX$. Thus, the conclusions from Theorems \ref{thm:SD_MBP}-\ref{thm:s_monotone} apply to every iteration of Procedure \ref{procedure:SD}.
\end{observe}

Having established that the repeated application of Survey Descent is locally well-defined, we now deduce the local Q-linear convergence induced by this procedure.

\begin{thm}
\textbf{(Q-Linear Convergence of Survey Points)\label{thm:q_linear}} Procedure \ref{procedure:SD} satisfies the following property:\begin{equation}
\underset{t\to\infty}{\limsup}\frac{\left\| \S^{t+1} - \bX \right\| ^{2}}{\left\| \S^{t}  - \bX \right\| ^{2}}\leq\left(1-\frac{\blam_{\min}\delta}{L}\right),\label{eq:s_result}
\end{equation}
where $\blam_{\min}$ is defined as in \eqref{eq:min_lam}.
\end{thm}
\begin{proof}
By Theorem \ref{thm:s_monotone} and Observation \ref{obs:allT}, there exists $K>0$ such that
\begin{equation*}
\left\| \S^{t+1} - \bX \right\|^{2}\leq\left(1-\frac{\blam_{\min}\delta}{L}\right)\left\| \S^{t} - \bX \right\|^{2}+K\left\| \S^{t} - \bX \right\|^{3} \ \forall \ t\in \nonneg.
\end{equation*}
Corollary \ref{cor:viable} implies that $\S^t$ is a valid survey. By Definition \ref{def:validS}, $\bX$ is not a valid survey because $f_i(\bx){=}f_j(\bx)$ for all $i,j$. Thus, $\S^t {\neq}\bX$, so $\left\| \S^{t} - \bX \right\|^2$ is nonzero and we can divide both sides by it.  We deduce
\begin{equation*}
\frac{\left\| \S^{t+1} - \bX \right\|^{2}}{\left\| \S^{t} - \bX \right\| ^{2}}\leq\left(1-\frac{\blam_{\min}\delta}{L}\right)+K\left\| \S^{t} - \bX\right\|  \ \forall \ t\in \nonneg.
\end{equation*}
\sloppy \edit{For $\S^{0}$ sufficiently close to $\bX$, this inequality guarantees that $\left\| \S^{t}{-}\bX\right\|$ monotonically decreases with $t$ and, in particular, $\left\| \S^{t}{-}\bX\right\|{<}\frac{\blam_{\min}\delta}{LK}$ for all $t\in\nonneg$.} Thus, $\frac{\left\| \S^{t+1} - \bX \right\|}{\left\| \S^{t} - \bX\right\|}{<}1$ implying $\underset{t\to\infty}{\lim}\left\| \S^{t}{-}\bX\right\|{=}0$. Taking the lim-sup with $t{\to}\infty$ on both sides of the above inequality then gives the desired result.
\end{proof}

Combining Theorem \ref{thm:q_linear} with the $L$-smoothness
assumptions (Definition \ref{def:strong_max}), we can immediately establish a loose
$R$-linear convergence result on the objective values.

\begin{cor}
\textbf{(R-Linear Convergence of Objective Values, Weak Version)\label{cor:Rlinear_weak}
} Consider Procedure \ref{procedure:SD}. Then, for any $K>\max_{j}\left\| \nabla f_{j}\left(\bx\right)\right\| _{2}$,
\[
f\left(s_{i}^{t}\right) - f(\bx) \leq K\left\| \S^{t} - \bX \right\| \ \forall\ i=1,\dots,k \text{, and } t\in\nonneg.
\]
Consequently, survey objective values $\left\{f\left(s_{i}^{t}\right)\right\}_{t\in\nonneg}$ converge R-linearly to $f(\bx)$ for all $i$.
\end{cor}
\begin{proof} Follows from the $L$-smoothness of the $f_i$-components and a routine application of the Cauchy-Schwarz inequality.
R-linear convergence follows from the Q-linear convergence of $\left\{\left\| \S^{t}{-}\bX \right\|\right\}_{t\in\nonneg}$ (Theorem \ref{thm:q_linear}).
\end{proof}

Compared to classical GD convergence guarantees on smooth objectives \citep[Theorem 10.29]{beck2017}, Corollary \ref{cor:Rlinear_weak} is rather weak: Firstly, the upper-bounding sequence is $\left\{ \left\| S^{t}{-}\bX\right\| \right\} _{t\in\nonneg}$ rather than the tighter $\left\{ \left\| S^{t}{-}\bX\right\|^{2}\right\}_{t\in\nonneg}$ seen in the smooth GD case; secondly, the constant $K$ depends on $\max_{j}\left\| \nabla f_{j}\left(\bx\right)\right\|_{2}$, which does not appear in the smooth GD case. A stronger function-value convergence guarantee for Survey Descent---more structurally similar to the aforementioned GD results---is indeed possible.

\begin{restatable}{thm}{Rlin}\label{thm:Rlin}

\textbf{(R-Linear Convergence of Objective Values, Strong Version)} Consider Procedure \ref{procedure:SD}. Fix $\bar{K}>\frac{5L}{\blam_{\min}}\left(1-\frac{\blam_{\min}\delta}{L}\right)$. Then,
\[
f\left(s_{i}^{t+1}\right) - f(\bx)\leq\bar{K}\left\| \S^{t} - \bX \right\| ^{2}\ \forall\ i=1,\dots,k \text{, and } t\in\nonneg.
\]
Consequently, survey objective values $\left\{f\left(s_{i}^{t}\right)\right\}_{t\in\nonneg}$ converge R-linearly to $f(\bx)$ for all $i$.
\end{restatable}
\begin{proof}
Follows from combining Theorem \ref{thm:approx_s}, Theorem \ref{thm:PR}, and the $L$-smoothness and $\delta$-strong convexity assumptions on the components $\left\{f_i\right\}_{i=1}^k$. We defer a detailed description to Appendix \ref{sec:strong_lconv}.
\end{proof}

Note that the rates of linear convergence presented in Theorems \ref{thm:q_linear} and \ref{thm:Rlin} are conservative.  Experiments suggest that these rate bounds are far from tight.

\section{Implementing Survey Descent}\label{sec:implement}
Stepping back from our theoretical analysis on strong $\C^2$ max functions $f$, we conclude with two remarks on the potential implementation of Survey Descent on general nonsmooth objectives $h$  (Definition \ref{def:SD_Method}).


\begin{rem}
\textbf{(Informal Heuristic for Survey Initialization)}\label{rem:init} We sketch an informal idea. First, run some initializing iterations of a standard method such as BFGS or a subgradient or proximal bundle method, and collect the gradients at these iterates. (As discussed, these iterates are typically points of differentiability.)  In practice, the iterates explore the nonsmooth landscape effectively.  In particular, on max-functions, they ``bounce'' between the active functions near the minimizer:  for BFGS, see our Figure \ref{fig:phenom} and \citet[Fig.\ 5.5,6]{lewis2009}, and for the subgradient method, see the figures and discussion on ``oscillations'' in \cite{bolte2020long}. As discussed, the convex hulls of nearby iterate gradients thus approximate the subdifferential at optimality, so we can estimate its dimension $d$ via singular value decomposition.  An initialization heuristic could then select a survey of $k=d+1$ iterates with robustly affinely independent gradients.  Implementing such heuristics is intricate:  we defer further discussion to future work.

\end{rem}

\begin{rem} \textbf{(Possible Acceleration of Survey Descent)} First, note that we can solve the $k$ constituent subproblems of Survey Descent (Definition \ref{def:SD_Method}) in parallel. Second, when all subproblems are feasible, we found in exploratory experiments that the following adjustment expedites Survey Descent's empirical convergence:
\begin{quote}
\textit{For all $i$, only update $s_i \leftarrow s_i^{+}$ when $h(s_i^{+}) < h(s_i)$; otherwise, keep $s_i$ the same in the subsequent survey.}
\end{quote}
This also tends to make Survey Descent compatible with a wider variety of initializing heuristics. For simplicity, we omit this additional enhancement from Definition \ref{def:SD_Method} and do not apply it in this paper's experiments.
\end{rem}

\subsection*{Acknowledgements}
Thanks are due to the Associate Editor and three anonymous referees for many helpful suggestions that significantly improved the initial manuscript.

\bibliographystyle{abbrvnat}
\bibliography{main}


\appendix

\part*{Appendix}
\renewcommand{\NotInAppendix}[1]{}

\section{Stronger R-Linear Convergence of Function Values (Theorem \ref{thm:Rlin})}\label{sec:strong_lconv}

In Section \ref{sec:lconv}, we analyzed the distance of the Survey Descent updates, $\S^+{=}\{s_i^+\}_{i=1}^k$, from the objective minimizer by examining the set of reference points $\{\ts_i\}_{i=1}^k$---which coincides with $\S^+$ up to $\mathrm{O}\left(\|\S\|^2\right)$ terms and is defined in \eqref{eq:s_tilde}. This appendix will prove Theorem \ref{thm:Rlin}---showing R-linear function value convergence---through an analogous strategy. To simplify the derivations, we adopt the following notation:\begin{equation}
\os_{i} =\proj\left[s_{i}-\frac{1}{L}\nabla f_{i}\left(s_{i}\right)\right] \text{ and }
\hs_{i} =\proj\left[s_{i}-\frac{\blam_{i}}{L}\nabla f_{i}\left(s_{i}\right)\right] \ \forall \ i=1,\dots,k. \label{eq:oracle_sol}
\end{equation}
In this notation, equation \eqref{eq:s_tilde} simplifies to
\begin{equation}
\ts_{i}\equiv\frac{1}{2-\blam_{i}}\hs_{i}+\sum_{j\neq i}\frac{\blam_{j}}{2-\blam_{i}}\os_{j}. \label{eq:s_tilde2}
\end{equation}
For any fixed $i$, we will prove Theorem \ref{thm:Rlin} in three steps:
\begin{enumerate}
\item \sloppy \textbf{Lemma \ref{lem:fi_to_f}:} Bounding the deviation between $f_{i}(\os_i)$  and $f(\os_i)$ because $\os_{i}$ does \textit{not} necessarily satisfy $f_{i}\left(\os_{i}\right){=} f\left(\os_{i}\right)$. Similarly, we also bound the difference between $f_{i}(\hs_i)$ and $f(\hs_i)$.
\item \textbf{Lemma \ref{lem:f_tilde}:} Using convexity and Theorem \ref{thm:PR} to bound $f\left(\ts_{i}\right)$.
\item \textbf{Lemma \ref{lem:f_tilde_diff}:} Bounding the difference between $f\left(s_{i}^{+}\right)$ and $f\left(\ts_{i}\right)$.
\end{enumerate}
We will again need to invoke the Taylor expansion \eqref{eq:taylor} to construct our bounds. However, we require a stronger version of Taylor's Theorem than that in Lemma \ref{lem:r_grad}---using ``Little-Oh'' rather than ``Big-Oh''. In particular, for a mapping $g:\mathbf{E}{\to}\F$ between two Euclidean spaces and letting $\left\| \cdot\right\|^p$ denote an arbitrary Euclidean norm raised to the $p$-th power, we use the notation $g\left(u\right)=\mathrm{o}\left(\left\| u\right\| ^{p}\right)$ to indicate the property of $g$ that---given \textit{any} constant $K{>}0$---there exists $U_{K}{>}0$ such that $\left|g\left(u\right)\right|{\leq}K\left\| u\right\| ^{p}$ for all $\left\| u\right\|{\leq}U_{K}$. By itself, we let $\mathrm{o}\left(\left\| u\right\| ^{p}\right)$ denote an element of the class of all functions with this property. Then, an identical argument as that in Lemma \ref{lem:r_grad} shows that the residual function in \eqref{eq:taylor} satisfies\begin{equation}
    r_i(x)=\mathrm{o}\left(\|x\|^2\right) \text{ and } \nabla r_{i}\left(x\right)=\mathrm{o}\left(\left\| x\right\|\right).\label{eq:r_grad_strong}
\end{equation}
We now prove our main result.
\begin{lem}
\label{lem:fi_to_f} Assume Setting \ref{setting:localSD}, and suppose $\bx{=}0$. Then, the following inequalities hold:
\begin{align*}
\left|f_{i}\left(\os_{i}\right)-f\left(\os_{i}\right) \right| &\leq \frac{L}{2}\left(1-\frac{\delta}{L}\right)\left\| s_{i} \right\| _{2}^{2}+\mathrm{o}\left(\left\| s_{i} \right\| _{2}^{2}\right),\\
\left|f_{i}\left(\hs_{i}\right)-f\left(\hs_{i}\right) \right|&\leq \frac{L}{2\blam_{i}}\left(1-\frac{\blam_{i}\delta}{L}\right)\left\| s_{i} \right\| _{2}^{2}+\mathrm{o}\left(\left\| s_{i} \right\| _{2}^{2}\right).
\end{align*}
\end{lem}
\begin{proof}
Using the Taylor expansion \eqref{eq:taylor}, observe for all pairs $i,j$ that
\begin{align*}
\left|f_{i}\left(\os_{i}\right)-f_{j}\left(\os_{i}\right)\right|
 &= \left|(a_i - a_j)^T \os_{i} + \frac{1}{2} \os_{i}^{T}\left(A_{i}-A_{j}\right)\os_{i}+r_{i}\left(\os_{i}\right)-r_{j}\left(\os_{i}\right)\right|\\
 & = \left|\frac{1}{2}\os_{i}^{T}\left(A_{i}-A_{j}\right)\os_{i}+r_{i}\left(\os_{i}\right)-r_{j}\left(\os_{i}\right)\right| \text{ (By $\os_i \in \U$ and \eqref{eq:Udef1})}\\
 & \leq\frac{L}{2}\left\| \os_{i}\right\| _{2}^{2}+\mathrm{o}\left(\left\| \os_{i}\right\| _{2}^{2}\right) \text{ (By $L$-smoothness and \eqref{eq:r_grad_strong})}\\
 & \leq\frac{L}{2}\left(1-\frac{\delta}{L}\right)\left\| s_{i}\right\| _{2}^{2}+\mathrm{o}\left(\left\| s_{i}\right\| _{2}^{2}\right)\ \text{ (Theorem \ref{thm:PR})}.
\end{align*}
Observe the right-hand side is independent of $j$, so the above still holds after replacing $f_j$ with $f{=}\max_j f_j$, which gives our desired result. An analogous argument for $\hs_{i}$\edit{, which replaces $\frac{1}{L}$ with $\frac{\blam_{i}}{L}$,} shows
\begin{equation*}
\left|f_{i}\left(\hs_{i}\right)-f\left(\hs_{i}\right)\right|\leq\frac{L}{2\blam_{i}}\left(1-\frac{\blam_{i}\delta}{L}\right)\left\| s_{i}\right\| _{2}^{2}+\mathrm{o}\left(\left\| s_{i}\right\| _{2}^{2}\right).
\end{equation*}
This completes the proof.
\end{proof}

\begin{lem}
\label{lem:f_tilde} Assume Setting \ref{setting:localSD}, and suppose $\bx{=}0$ and $f(\bx){=}0$. Then, the following inequality holds:\begin{equation}
f\left(\ts_{i}\right) \leq\frac{L}{\blam_{\min}}\left(1-\frac{\delta\blam_{\min}}{L}\right)\max_{j=1,\dots,k}\left\| s_{j} \right\| _{2}^{2}+\mathrm{o}\left(\left\| \S \right\| ^{2}\right),\label{eq:f_aux}
\end{equation}
where $\ts_i$ and $\blam_{\min}$ are as defined in equations \eqref{eq:s_tilde} and \eqref{eq:min_lam} respectively.
\end{lem}
\begin{proof}
Using \eqref{eq:s_tilde2}, Theorem \ref{thm:approx_s}, and convexity, we deduce
\begin{equation}
f\left(\ts_{i}\right) =  f\left(\frac{1}{2-\blam_{i}}\hs_{i}+\sum_{j\neq i}\frac{\blam_{j}}{2-\blam_{i}}\os_{j}\right) \leq  \frac{1}{2-\blam_{i}}f\left(\hs_{i}\right)+\sum_{j\neq i}\frac{\blam_{j}}{2-\blam_{i}}f\left(\os_{j}\right). \label{eq:convex_up}
\end{equation}
Define the quantity
\begin{equation*}
\Delta\equiv\frac{L}{2\blam_{\min}}\left(1-\frac{\blam_{\min}\delta}{L}\right)\left(\frac{1}{2-\blam_{i}}\left\| s_{i}\right\| _{2}^{2}+\sum_{j\neq i}\frac{\blam_{j}}{2-\blam_{i}}\left\| s_{j}\right\|_{2}^{2}\right).
\end{equation*}
Applying Lemma \ref{lem:fi_to_f} to \eqref{eq:convex_up} \edit{combined with the facts that $\blam_{\min} \leq \blam_{j} \leq 1$ for all $j$} gives\begin{equation}
f\left(\ts_{i}\right) \leq \underbrace{ \frac{1}{2-\blam_{i}}f_{i}\left(\hs_{i}\right)+\sum_{j\neq i}\frac{\blam_{j}}{2-\blam_{i}}f_{j}\left(\os_{j}\right)}_{\left(\ast\right)} + \Delta + \mathrm{o}\left(\left\| \S\right\| ^{2}\right).\label{eq:convex_up2}
\end{equation}
Theorem \ref{thm:PR} implies the following bound on $\left(\ast\right)$:
\begin{align*}
\left(\ast\right) & \leq\frac{1}{2-\blam_{i}}\frac{L}{2\blam_{i}}\left(1-\frac{\delta\blam_{i}}{L}\right)\left\| s_{i}\right\| ^{2}+\frac{L}{2}\left(1-\frac{\delta}{L}\right)\sum_{j\neq i}\frac{\blam_{j}}{2-\blam_{i}}\left\| s_{j}\right\| _{2}^{2}\leq\Delta.
\end{align*}
Substituting the above into \eqref{eq:convex_up2} gives the desired
\begin{align*}
f\left(\ts_{i}\right) & \leq2\Delta+\mathrm{o}\left(\left\| \S\right\| _{2}^{2}\right) =\frac{L}{\blam_{\min}}\left(1-\frac{\delta\blam_{\min}}{L}\right)\left[\frac{1}{2-\blam_{i}}\left\| s_{i}\right\| ^{2}+\sum_{j\neq i}\frac{\blam_{j}}{2-\blam_{i}}\left\| s_{j}\right\| _{2}^{2}\right]+\mathrm{o}\left(\left\| \S\right\| ^{2}\right)\\
 & \leq\frac{L}{\blam_{\min}}\left(1-\frac{\delta\blam_{\min}}{L}\right)\max_{j=1,....,k}\left\| s_{j}\right\| _{2}^{2}+\mathrm{o}\left(\left\| \S\right\| ^{2}\right),
\end{align*}
which is our desired result.
\end{proof}

\begin{lem} \label{lem:f_tilde_diff} Assume Setting \ref{setting:localSD}, and suppose $\bx{=}0$ and $f(\bx){=}0$. Then, the following inequality holds:
\begin{equation*}
f\left(s_{i}^{+}\right) - f\left(\ts_{i}\right) \leq \frac{4L}{\blam_{\min}}\left(1-\frac{\blam_{\min}\delta}{L}\right)\max_{j=1,...,k}\left\| s_{i}\right\| _{2}^{2}+\mathrm{O}\left(\left\| \S\right\| ^{3}\right),
\end{equation*}
where $\ts_i$ and $\blam_{\min}$ are as defined in equations \eqref{eq:s_tilde} and \eqref{eq:min_lam} respectively.
\end{lem}
\begin{proof}
Consider any fixed $i$. By Theorem \ref{thm:SD_MBP}, all constraints of $\left(P_{i}^{\S}\right)$
are tight, so $s_{i}^{+}$ satisfies the constraints of Observation \ref{obs:SD_MOS} with equality for all $j\neq i$:\begin{equation}
\ell_{s_{j}}^{f_{j}}\left(s_{i}^{+}\right)+\frac{L}{2}\left\| s_{i}^{+}-s_{j}\right\| _{2}^{2}=\ell_{s_{i}}^{f_i}\left(s_{i}^{+}\right).\label{eq:lem51_lin}
\end{equation}
Define $\xi_{i} \equiv s_{i}^{+}{-}\ts_{i}$. In other words, $\xi_{i}$ is the $\mathrm{O}\left(\|\S\|^2\right)$ term in equation \eqref{eq:SD_iter}.
Consider any fixed $j$. Substituting $s_i^+ = \ts_i + \xi_i$ into \eqref{eq:lem51_lin} and expanding the quadratic term gives
\begin{equation*}
\ell_{s_{j}}^{f_{j}}\left(\ts_{i}\right)+
\nabla f_j(s_j)^{T}\xi_{i}+\frac{L}{2}\left\| \ts_{i}-s_{j}\right\| _{2}^{2}+L\left\langle \ts_{i}-s_{j},\xi_{i}\right\rangle +\frac{L}{2}\left\| \xi_{i}\right\|^{2}
=\ell_{s_{i}}^{f_{i}}\left(\ts_{i}\right)+\nabla f_i(s_i)^{T}\xi_{i}.
\end{equation*}
Applying the Taylor expansion \eqref{eq:taylor} to $\nabla f_i(s_i)$ and $\nabla f_j(s_j)$ gives
\begin{multline*}
\ell_{s_{j}}^{f_{j}}\left(\ts_{i}\right)+
\left(a_j + A_j s_j + \nabla r_j(s_j)\right)^T\xi_{i}+\frac{L}{2}\left\| \ts_{i}-s_{j}\right\| _{2}^{2}+L\left\langle \ts_{i}-s_{j},\xi_{i}\right\rangle +\frac{L}{2}\left\| \xi_{i}\right\|^{2}
\\
=\ell_{s_{i}}^{f_{i}}\left(\ts_{i}\right)+\left(a_i + A_i s_i + \nabla r_i(s_i)\right)^{T}\xi_{i}.
\end{multline*}
Note that the terms $\nabla r_{j}(s_j)$, $\nabla r_i(s_i)$, and $\ts_{i}$ are all $\mathrm{O}\left(\left\| \S\right\| \right)$ by Lemma \ref{lem:r_grad} and \eqref{eq:s_tilde} and \edit{so is $s_i$ by definition}. Thus, collecting $\mathrm{O}\left(\left\| \S\right\| ^{3}\right)$ terms simplifies the above to
\begin{equation}
\ell_{s_{j}}^{f_{j}}\left(\ts_{i}\right)+a_{j}^{T}\xi_{i}+\frac{L}{2}\left\| \ts_{i}-s_{j}\right\| _{2}^{2}+\mathrm{O}\left(\left\| \S\right\| ^{3}\right)=\ell_{s_{i}}^{f_{i}}\left(\ts_{i}\right)+a_{i}^{T}\xi_{i}. \label{eq:lem51_simp}
\end{equation}
Multiplying \eqref{eq:lem51_simp} by $\blam_{j}$, summing over $j\neq i$, and using identities
\eqref{eq:lam_sum1} and \eqref{eq:lam_sum0} gives
\begin{align*}
\sum_{j\neq i}\blam_{j}\ell_{s_{j}}^{f_{j}}\left(\ts_{i}\right)-\blam_{i}a_{i}^{T}\xi_{i}+\frac{L}{2}\sum_{j\neq i}\left(\blam_{j}\left\| \ts_{i}-s_{j}\right\| _{2}^{2}\right)+\mathrm{O}\left(\left\| \S\right\| ^{3}\right) & =\left(1-\blam_{i}\right)\ell_{s_{i}}^{f_{i}}\left(\ts_{i}\right)+\left(1-\blam_{i}\right)a_{i}^{T}\xi_{i}.
\end{align*}
Adding $\blam_i \ell^{i}_{s_i}(\ts_i)$ to both sides \edit{and again applying Lemma \ref{lem:r_grad}} leads to
\begin{align*}
\sum_{j=1}^{k}\blam_{j}\ell_{s_{j}}^{f_{j}}\left(\ts_{i}\right)+\frac{L}{2}\sum_{j\neq i}\left(\blam_{j}\left\| \ts_{i}-s_{j}\right\| _{2}^{2}\right)+\mathrm{O}\left(\left\| \S\right\| ^{3}\right) & =\ell_{s_{i}}^{f_{i}}\left(\ts_{i}\right)+a_{i}^{T}\xi_{i} =\ell_{s_{i}}^{f_{i}}\left(s_{i}^{+}\right)\underbrace{-\left(A_{i}s_{i}+\nabla r_{i}\left(s_{i}\right)\right)^{T}\xi_{i}}_{\mathrm{O}\left(\left\| \S\right\| ^{3}\right)}.
\end{align*}
Again collecting  $\mathrm{O}\left(\left\| \S\right\| ^{3}\right)$ terms and adding a non-negative $\frac{\blam_{i} L}{2} \|\ts_i - s_i\|_2^2$ term to the left-hand side leads to the following inequality:
\begin{align}
\underbrace{\sum_{j=1}^{k}\blam_{j}\ell_{s_{j}}^{f_{j}}\left(\ts_{i}\right)}_{\left(\sun\right)}+\frac{L}{2}\sum_{j=1}^{k}\left(\blam_{j}\left\| \ts_{i}-s_{j}\right\| _{2}^{2}\right)+\mathrm{O}\left(\left\| \S\right\| ^{3}\right) & \geq\ell_{s_{i}}^{f_{i}}\left(s_{i}^{+}\right)\label{eq:f_lin_bound}
\end{align}
For $\left(\sun\right)$, the assumed $\delta$-strong convexity on $f_j$ (Definition \ref{def:strong_max}) implies
\begin{equation}
\sum_{j=1}^{k}\blam_{j}\ell_{s_{j}}^{f_{j}}\left(\ts_{i}\right)\leq\sum_{j=1}^{k}\blam_{j}\left(f_{j}\left(\ts_{i}\right)-\frac{\delta}{2}\left\| \ts_{i}-s_{j}\right\| _{2}^{2}\right)\leq\underbrace{\max_{j=1,\dots,k}f_{j}\left(\ts_{i}\right)}_{=f\left(\ts_{i}\right)}-\frac{\delta}{2}\sum_{j=1}^{k}\blam_{j}\left\| \ts_{i}-s_{j}\right\| _{2}^{2},\label{eq:f_lin_bound2}
\end{equation}
where the right-most inequality follows from \eqref{eq:lam_sum1}.
Substituting \eqref{eq:f_lin_bound2} into \eqref{eq:f_lin_bound} leads to
\begin{equation}
    \begin{aligned}
    f\left(\ts_{i}\right)+\frac{L}{2}\left(1-\frac{\delta}{L}\right)\sum_{j=1}^{k}\left(\blam_{j}\left\| \ts_{i}-s_{j}\right\| _{2}^{2}\right)+\mathrm{O}\left(\left\| \S\right\| ^{3}\right) & \geq\ell_{s_{i}}^{f_{i}}\left(s_{i}^{+}\right) \geq f_{i}\left(s_{i}^{+}\right)-\frac{L}{2}\left\| s_{i}^{+}-s_{i}\right\| _{2}^{2},
    \end{aligned}\label{eq:lem51_f}
\end{equation}
where the second inequality follows from the $L$-smoothness of $f_i$ (Definition \ref{def:strong_max}). By Theorem \ref{thm:strictval}, $s_i^+$ is a valid survey point, so $f_{i}\left(s_{i}^{+}\right)=f\left(s_{i}^{+}\right)$. Substituting this into \eqref{eq:lem51_f} and rearranging gives
\begin{equation}
f\left(s_{i}^{+}\right)\leq f\left(\ts_{i}\right)+\frac{L}{2}\left(1-\frac{\delta}{L}\right)\underbrace{\sum_{j=1}^{k}\left(\blam_{j}\left\| \ts_{i}-s_{j}\right\| _{2}^{2}\right)}_{\left(\dagger\right)}+\frac{L}{2}\underbrace{\left\| s_{i}^{+}-s_{i}\right\| _{2}^{2}}_{\left(\star\right)}+\mathrm{O}\left(\left\| \S\right\| ^{3}\right)\label{eq:f_lin_bound3}
\end{equation}
We use Jensen's inequality to bound $\left(\star\right)$:
\begin{align*}
\left(\star\right)\leq\max_{j=1,...,k}\left\| s_{i}^{+}-s_{j}\right\| _{2}^{2} & =\max_{j=1,...,k}\left\| \frac{1}{2}\left(2s_{i}^{+}\right)+\frac{1}{2}\left(-2s_{j}\right)\right\| _{2}^{2}\leq\frac{1}{2}\left\| \left(2s_{i}^{+}\right)\right\| _{2}^{2}+\frac{1}{2}\max_{j=1,...,k}\left\| \left(-2s_{j}\right)\right\| _{2}^{2}.
\end{align*}
Hence,
\begin{align*}
\left(\star\right) & \leq2\left(\left\| s_{i}^{+}\right\| _{2}^{2}+\max_{j=1,...,k}\left\| s_{j}\right\| _{2}^{2}\right)\\
 & =2\left(2-\frac{\blam_{\min}\delta}{L}\right)\max_{j=1,....,k}\left\| s_{j}\right\| _{2}^{2}+\mathrm{O}\left(\left\| \S\right\|^{3}\right)\text{ (Theorem \ref{thm:s_monotone})}\\
 & \leq4\max_{j=1,....,k}\left\| s_{j}\right\| _{2}^{2}+\mathrm{O}\left(\left\| \S\right\| ^{3}\right).
\end{align*}
An analogous argument (using Lemma \ref{lem:st_bound} on $\ts_i$ instead of Theorem \ref{thm:s_monotone} on $s_i^+$) gives
\begin{equation*}
\left(\dagger\right)\leq\max_{j=1,...,k}\left\| \ts_{i}-s_{j}\right\| _{2}^{2}\leq4\max_{j=1,....,k}\left\| s_{j}\right\| _{2}^{2}.
\end{equation*}
Substituting these bounds on $\left(\dagger\right)$ and $\left(\star\right)$
into \eqref{eq:f_lin_bound3}, we deduce
\begin{align*}
f\left(s_{i}^{+}\right) & \leq f\left(\ts_{i}\right)+2L\left(1-\frac{\delta}{L}\right)\max_{j=1,...,k}\left\| s_{i}\right\| _{2}^{2}+2L\max_{j=1,...,k}\left\| s_{i}\right\| _{2}^{2}+\mathrm{O}\left(\left\| \S\right\| ^{3}\right)\\
 & =f\left(\ts_{i}\right)+4L\left(1-\frac{\delta}{2L}\right)\max_{j=1,...,k}\left\| s_{i}\right\| _{2}^{2}+\mathrm{O}\left(\left\| \S\right\| ^{3}\right).
\end{align*}

If $k=1$, Survey Descent reduces to GD and the result is trivial from \citet[Theorem 10.29]{beck2017} or \citet[Theorem 2.1.15]{nesterov2003introductory}. When $k \geq 2$, the property \eqref{eq:lam_sum1} implies $\blam_{\min}\leq\frac{1}{k}\leq\frac{1}{2}$. Thus, $\left(1-\frac{\delta}{2L}\right){\leq}\left(1-\frac{\blam_{\min}\delta}{L}\right)$.
We then deduce
\begin{equation*}
f\left(s_{i}^{+}\right)\leq f\left(\ts_{i}\right)+\frac{4L}{\blam_{\min}}\left(1-\frac{\blam_{\min}\delta}{L}\right)\max_{j=1,...,k}\left\| s_{i}\right\| _{2}^{2}+\mathrm{O}\left(\left\| \S\right\| ^{3}\right).
\end{equation*}
This completes the proof. 
\end{proof}

\Rlin*
\begin{proof}
Without loss of generality, assume $\bx{=}0$ and $f(\bx){=}0$. Consider any fixed $i$ and $t$. For simplicity, adopt the notation $s_i{=}s_i^t$ and $s_i^+{=}s_i^{t+1}$. Adding the inequalities in Lemmas \ref{lem:f_tilde} and \ref{lem:f_tilde_diff} gives
\begin{align*}
f\left(s_{i}^{+}\right) & \leq\frac{5L}{\blam_{\min}}\left(1-\frac{\blam_{\min}\delta}{L}\right)\max_{j=1,...,k}\left\| s_{j}\right\| _{2}^{2}+\underbrace{\mathrm{o}\left(\left\| \S\right\| ^{2}\right)+\mathrm{O}\left(\left\| \S\right\| ^{3}\right)}_{=\mathrm{o}\left(\left\| \S\right\| ^{2}\right)}.
\end{align*}
Then, for any $K>0$ and for all $\S$ sufficiently close to $\bX$,
we have
\begin{align*}
f\left(s_{i}^{+}\right) & \leq\frac{5L}{\blam_{\min}}\left(1-\frac{\blam_{\min}\delta}{L}\right)\left\| \S\right\| ^{2}+K\left\| \S\right\| ^{2}=\left(\frac{5L}{\blam_{\min}}\left(1-\frac{\blam_{\min}\delta}{L}\right)+K\right)\left\| \S\right\| ^{2}.
\end{align*}
Denote $\bar{K}\equiv\frac{5L}{\blam_{\min}}\left(1-\frac{\blam_{\min}\delta}{L}\right)+K>0$.
Since $K$ can be any strictly positive constant, choosing $K$ is
equivalent to choosing an arbitrary $\bar{K}>\frac{5L}{\blam_{\min}}\left(1-\frac{\blam_{\min}\delta}{L}\right)$.
Under such a choice,
\begin{equation*}
f\left(s_{i}^{+}\right)\leq\bar{K}\left\| \S\right\| ^{2},
\end{equation*}
which is our desired result. The general case follows from applying the change of variables $x{\leftarrow}(x{-}\bx)$ and function-value shift $f(\cdot){\leftarrow}\left(f(\cdot){-}f(\bx)\right)$.
\end{proof}

\ifrevise \else 

    \section{Geometric Proof of Theorem \ref{thm:SD_MBP}}\label{sec:thm_geometry} 
    
    Although Theorem \ref{thm:SD_MBP} can be proven by verifying first-order optimality conditions and applying the implicit function theorem, we present here a more direct, computationally revealing approach that highlights the underlying geometry of Survey Descent. In particular, after change-of-variables and tightness-of-constraints arguments (presented later), Observation \ref{obs:SD_MOS} reduces to a simple, geometric \textit{multiple balls problem.}
    
    \begin{defn}
    \textbf{(Multiple Balls Problem, $\MBP\param$)\label{def:MBP}} \sloppy For some
    target vector $y\in\R^{n}$, ball centers $Z=\left[z_{1},\dots,z_{\bk}\right]\in\R^{n\times\bk}$
    where $\bk\leq n$, and non-negative radii $r=\left[r_{1},\dots,r_{\bk}\right]^{T}\in\R^{\bk}_{+}$,
    define the \textit{Multiple Balls Problem (MBP)} as 
    \begin{align}
    \min_{x\in\R^{n}}\ & \left\| x-y\right\| _{2}^{2} \nonumber\\
    \mathrm{s.t.}\ & \left\| x-z_{j}\right\| _{2}^{2} = r_{j}^{2}\ \forall\ j=1,\dots,\bk. \label{eq:MBP_cons}
    \end{align}
    Denote the problem as $\MBP\param$ and its \textit{feasible region} as $\B\spec$. We will refer to the tuple $\param$ as the \textit{parameters} and the sub-tuple $\spec$ as the \textit{ball-specifications} of the MBP.
    \end{defn}
    Thus, MBP simply projects $y$ onto the boundary-intersection, $\B\spec$, of $\bk$ balls of radii $\left\{ r_{j}\right\} _{j=1}^{\bk}$---each centered at the corresponding $\left\{ z_{j}\right\} _{j=1}^{\bk}$.
    This is illustrated in Figure \ref{fig:MBP_params}. 
    
    Since $y$, $Z$, and $r$ all reside in Euclidean spaces, we can also consider $\spec$ and $\param$ elements of some product space. Let $\left\Vert \cdot\right\Vert $ denote an arbitrary Euclidean norm. This will allow us to formalize the notion of distance in of parameters and ball-specifications in their respective spaces. We are particularly interested in the behavior of $\MBP\param$ where $\param$ is near some set of \textit{base parameters}.
    
    \begin{defn}
    \textbf{(Base Parameters)\label{def:base}} We refer to a set of parameters
    as \textit{base parameters} and denote them as $\bparam$
    if they satisfy the following properties:
    \begin{itemize}
    \item The columns of $\bZ=\left[\bz_{1},\dots,\bz_{\bk}\right]$ are linearly
    independent.
    \item $\br_{j}^{2}=\left\| \bz_{j}\right\| _{2}^{2}$.
    \item There exist strictly-positive weights $\clam=\left[\clam_{1},\dots,\clam_{\bk}\right]^{T}\in\R_{++}^{\bk}$
    such that $\bar{y}=-\sum_{j=1}^{\bk}\clam_{j}\bar{z}_{j}.$
    \end{itemize}
    \end{defn}
    Note that $\clam$ is necessarily unique due to the linear independence of $\bZ$, and their strict-positivity are connected to the nondegenerate first-order optimality that we will later show for $\MBP\bparam$.
    
    The remainder of Appendix \ref{sec:thm_geometry} is organized as follows. For parameters $\param$ near some base parameters $\bparam$, Section \ref{sec:MBP_sol} derives an explicit form solution for $\MBP\param$, Section \ref{sec:MBP_ineq} shows that the solution is the same if one were to solve $\MBP\param$ with inequality rather than equality constraints, and Section \ref{sec:MBP_SD} connects the MBPs to Survey Descent.
    
    \subsection{Local, Explicit-Form Solution to the Multiple Balls Problem}\label{sec:MBP_sol}
    
    In Theorem \ref{thm:MBP_sol}, we will build the explicit-form solution for $\MBP\param$---for parameters near $\bparam$---from $\B\spec$ as well as the following geometric objects (whose well-definedness and other properties will be derived in Lemmas \ref{lem:A}-\ref{lem:B}):
    \begin{itemize}
        \item The \textit{central affine subspace}
            \begin{equation}
            \A\spec \equiv \left\{ x\in\R^{n}:\ 2z_{j}^{T}x-\left\| z_{j}\right\| _{2}^{2}+r_{j}^{2}\text{ equal  for all } j=1,\dots,\bk\right\}. \label{eq:CAS}
            \end{equation}
        \item The \textit{core}
            \begin{equation}
            c\spec\equiv\projA{\spec}\left(z_{j}\right)\ \text{ for any }\ j=1,\dots,\bk, \label{eq:core}
            \end{equation}
            where $\projA{\spec}$ denotes the projection onto $\A\spec$.
        \item The \textit{margin}
            \begin{equation}
            \rho\spec\equiv\sqrt{r_{j}^{2}-\left\| z_{j}-c\spec\right\| _{2}^{2}} \ \text{ for any }\ j=1,\dots,\bk. \label{eq:margin}
            \end{equation}
    \end{itemize}
    They are collectively visualized in Figure \ref{fig:MBP_objects}.
    
    \begin{lem}\textbf{(Local Properties of the Central Affine Subspace)}\label{lem:A}
    For any ball-specifications $\spec$ sufficiently close to base ball-specifications $\bspec$, the following properties hold:
    \begin{enumerate}
        \item $\A\spec$ is non-empty with dimension $n{-}\bk{+}1$.
        \item For all $i\neq j$, the direction $(z_i - z_j)$ is orthogonal to the linear subspace parallel to $\A\spec$.
        \item For $w\in\R^n$, the projection $\projA{\spec}(w)$ is well-defined and varies smoothly with $\spec$ and $w$.
    \end{enumerate}
    \end{lem}
    \begin{proof}
    Observe that we can rewrite \eqref{eq:CAS} as the following $\left(\bk{-}1\right){\times}n$ linear system:\begin{equation}
    \A\spec=\left\{ x\in\R^{n}:\ M\spec x=b\spec \right\} ,\label{eq:CAS_sys}
    \end{equation}
    where 
    \[
    M\spec\equiv2\begin{bmatrix}\left(z_{2}-z_{1}\right)^{T}\\
    \vdots\\
    \left(z_{\m}-z_{1}\right)^{T}
    \end{bmatrix}
    \text{ and }
    b\spec\equiv\begin{bmatrix}\left(\left\| z_{2}\right\| _{2}^{2}-r_{2}^{2}\right)-\left\| z_{1}\right\| _{2}^{2}+r_{1}^{2}\\
    \vdots\\
    \left(\left\| z_{\bk}\right\| _{2}^{2}-r_{\bk}^{2}\right)-\left\| z_{1}\right\| _{2}^{2}+r_{1}^{2}
    \end{bmatrix}.\
    \]
    
    Consider $\A\bspec$. Since $\left\{ \bz_{j}\right\} _{j=1}^{\bk}$ are linearly independent,
    $\left\{ \bz_{j}-\bz_{1}\right\} _{j=2}^{\bk}$---which are the rows of $M\bspec$---are also linearly independent. For $\spec$ sufficiently close to $\bspec$, the rows of $M\spec$ will remain linearly independent. In this case, the linear system \eqref{eq:CAS_sys} is then necessarily solvable: $\A\spec$ is non-empty. Let $x_0$ be such a solution. Then, $\A\spec$ consists of all vectors $x \in \R^n$ such that $x{-}x_0$ is orthogonal to the $m{-}1$ linearly independent rows of $M\spec$. Thus, $\Dim\left\{\A\spec\right\}{=}n{-}\bk{+}1$. 
    
    The vectors $\left\{(z_i - z_j)\right\}_{j\neq i}$ clearly all lie in the rowspace of $M\spec$. Hence, the vectors$\left\{(z_i - z_j)\right\}_{j\neq i}$ are orthogonal to all direction in $\A\spec$.
    
    Finally, write $M$ instead of $M\spec$ for simplicity. Since $M$ has independent rows, the projection of $w$ onto $\A\spec$ is explicitly given by
    \begin{equation*}
    \projA{\spec}\left(w\right)=\left(I-M\left(M^{T}M\right)^{-1}M^{T}\right)w+M\left(M^{T}M\right)^{-1}b,
    \end{equation*}
    which varies smoothly with $\spec$ and $w$.
    \end{proof}
    \begin{figure}[t]
    \subfloat[\label{fig:MBP_params}]{\includegraphics[width=3.25in]{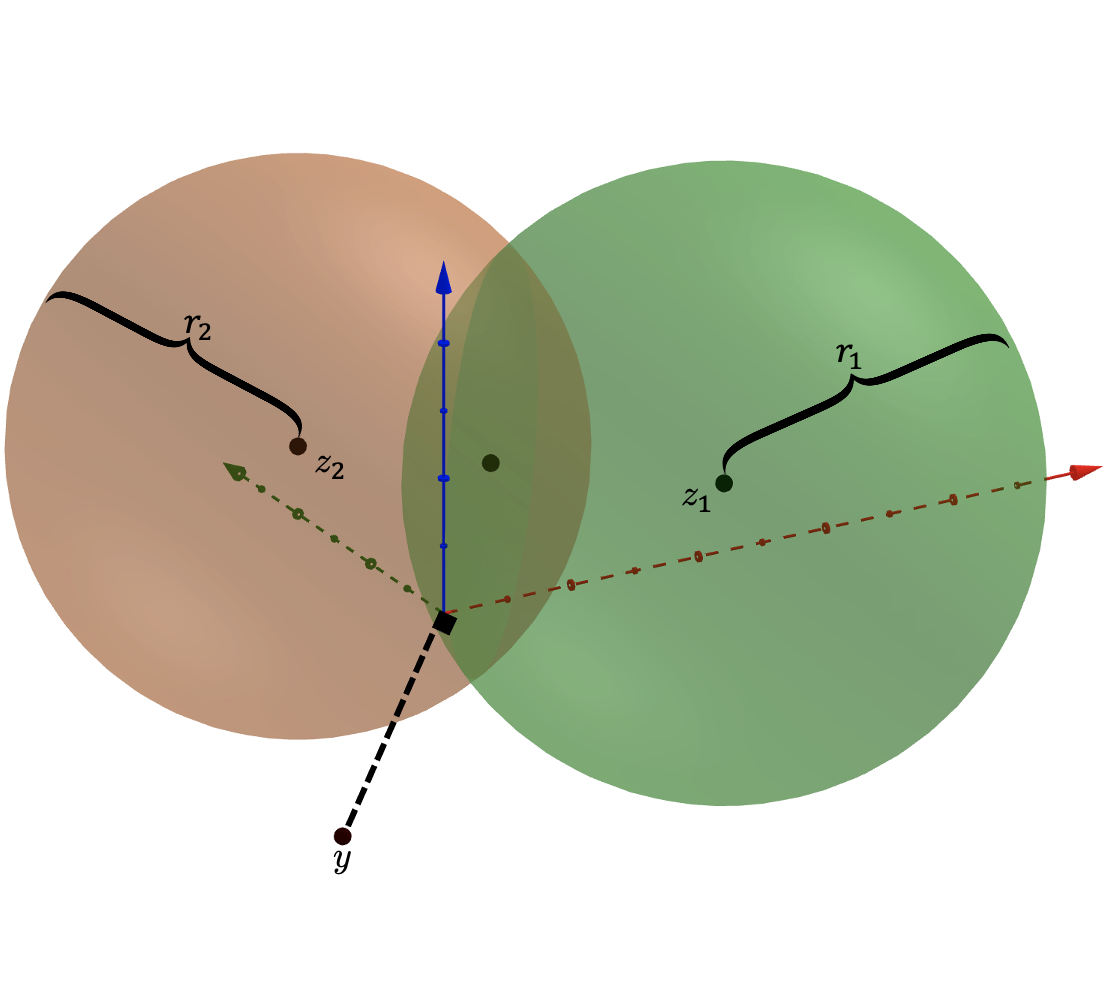}}\hfill{}\subfloat[\label{fig:MBP_objects}]{\includegraphics[width=3.25in]{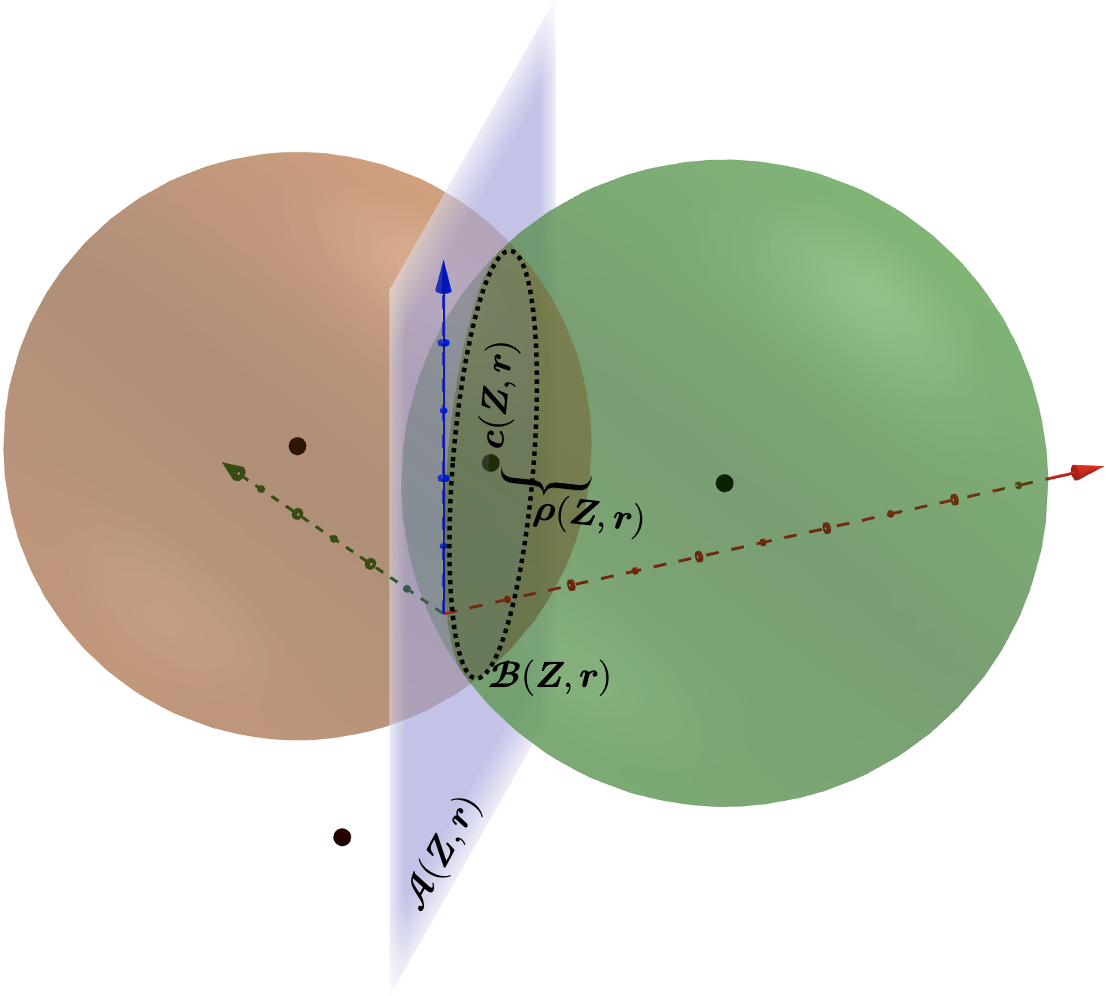}}
    
    \caption{\textbf{Visualization of the Multiple Balls Problem.} Given a target $y\in\protect\R^{n}$, ball centers $Z=\left[z_{1},\dots,z_{\protect\bk}\right]\in\protect\R^{n\times\protect\bk}$, and corresponding radii $r=\left[r_{1},\dots,r_{\protect\bk}\right]^{T}\in\protect\R^{\protect\bk}$, the above displays the geometries explored in Appendix \ref{sec:thm_geometry} for $n=3$ and $\protect\bk=2$. Panel \textbf{(a)} shows the Multiple Balls Problem in which the point $y$ is projected onto the boundary-intersection of the $\protect\bk$-balls specified by $\spec$. Panel \textbf{ (b)} illustrates the four fundamental structures explored in Appendix \ref{sec:thm_geometry}: the central affine subspace $\protect\A\spec$, the core $c\spec$, the margin $\rho\spec$, 
    and the feasible region $\protect\B\spec$.
    }
    \end{figure}
    
    \FloatBarrier
    
    \begin{lem}
    \textbf{\label{lem:core} (Local Properties of Core)} For any ball-specifications $\spec$ sufficiently close to base ball-specifications $\bspec$, the core in \eqref{eq:core} is well-defined and varies smoothly with $\spec$.
    \end{lem}
    \begin{proof} By the second fact in Lemma \ref{lem:A}, for any pair $i\neq j$, the difference $(z_i{-}z_j)$ is orthogonal to the linear subspace parallel to $\A\param$. Therefore, $\projA{\param}(z_i)$ and $\projA{\param}(z_j)$ are the same. Hence, \eqref{eq:core} is well-defined. The smoothness of the $c\spec$ with respect to $\spec$ then follows from \eqref{eq:core} and the third fact of Lemma \ref{lem:A}.
    \end{proof}
    
    \begin{lem}
    \textbf{(The Base Case)} \label{lem:linA}
    At base ball-specifications, the central affine subspace is a linear subspace and the core is non-zero.
    \end{lem}
    \begin{proof}
    By Definition \ref{def:base}, $\br_{j}^{2}=\left\| \bz_{j}\right\| _{2}^{2}$
    for all $j$. Thus, \eqref{eq:CAS} simplifies to
    \begin{equation*}
    \A\bspec=\left\{ x\in\R^{n}:\ 2\bz_{j}^{T}x\text{ equal  for all } j=1,\dots,\bk\right\}.
    \end{equation*}
    Clearly, the origin satisfies the RHS conditions above, so $0\in\A\bspec$.
    
    \sloppy Next, assume for contradiction that $0{=}c\bspec$. By definition, $c\bspec$ is the projection of $\bz_j$ onto $\A\bspec$ for any $j$, so the vectors $\left\{ \bz_{j}-c\bspec\right\}_{j=1}^{\bk}{=}\left\{ \bz_{j}\right\} _{j=1}^{\bk}$ must be orthogonal to $\A\bspec$. The linear independence of $\left\{ \bz_{j}\right\}_{j=1}^{\bk}$ implies $\Dim\left\{ \A\bspec\right\}{\leq}n{-}\bk$, which contradicts the first fact of Lemma \ref{lem:A}. Thus, $0{\neq}c\bspec$.
    \end{proof}
    
    \begin{lem}
    \textbf{(Local Properties of the Margin)\label{lem:margin}} 
    For ball-specifications $\spec$ near some base ball-specifications $\bspec$, the margin in \eqref{eq:margin} is well-defined, strictly-positive, and varies smoothly with $\spec$.
    \end{lem}
    \begin{proof}
    It suffices to show that the square of the margin,\begin{equation}
    r_{j}^{2}-\left\| z_{j}-c\spec\right\| _{2}^{2}, \label{eq:margin_sq}
    \end{equation}
    does not depend on the index $j$, is strictly-positive, and varies smoothly with $\param$. First, consider the core and central affine subspace exactly at base ball-specifications $\bspec$. By Lemma
    \ref{lem:linA}, since $0 \in \A\bspec$ and $\projA{\bspec}\left(\bz_{j}\right){=}c\bspec{\neq}0$
    for all $j$, the definition of a projection gives
    \begin{equation*}
    \left\| \bz_{j}-c\bspec\right\| _{2}<\left\| \bz_{j}-0\right\| _{2}=\br_{j}\ \forall\ j=1,\dots,\bk.
    \end{equation*}
    Thus, \eqref{eq:margin_sq} is strictly positive, which implies the argument of the square-root in \eqref{eq:margin} is strictly positive. Next, observe that
    \begin{align*}
    r_{j}^{2}-\left\| z_{j}-c\spec\right\| _{2}^{2} & =-\left\| c\spec\right\| _{2}^{2}+\underbrace{2z_{j}^{T}c\spec-\left\| z_{j}\right\| _{2}^{2}+r_{j}^{2}}_{(\star)_j} \ \forall\  j=1,\dots,\bk.
    \end{align*}
    Since $c\spec{\in}\A\spec$, the term $(\star)_j$ is equal for all $j$. Therefore, \eqref{eq:margin_sq} is $j$-independent, and $\rho\spec$ in \eqref{eq:margin} is well-defined. The smoothness of $\rho\spec$  with respect to $\spec$ then follows directly from \eqref{eq:margin} and the smoothness of the core (Lemma \ref{lem:core}).
    \end{proof}
    
    \begin{lem}
    \label{lem:B}\textbf{(Local Properties of the Feasible Region)} For ball-specifications
    $\spec$ near some base ball-specifications $\bspec$, the feasible region is nonempty and satisfies\begin{equation}
    \B\spec=\left\{ x\in\A\spec:\ \left\| x-c\spec\right\| _{2}=\rho\spec\right\}.\label{eq:Bdef2} 
    \end{equation}
    \end{lem}
    \begin{proof}
    \noindent\textit{Equivalence of \eqref{eq:MBP_cons} and \eqref{eq:Bdef2}:}
    Expanding and rearranging \eqref{eq:MBP_cons} gives
    \begin{align*}
    \left\| x-z_{j}\right\| _{2}=r_{j}\ \forall\ j & \iff\left\| x\right\| _{2}^{2}=2z_{j}^{T}x-\left\| z_{j}\right\| _{2}^{2}+r_{j}^{2}\ \forall\ j\\
    & \implies2z_{j}^{T}x-\left\| z_{j}\right\| _{2}^{2}+r_{j}^{2}=2z_{i}^{T}x-\left\| z_{i}\right\| _{2}^{2}+r_{i}^{2}\ \forall\ i\neq j.
    \end{align*}
    Hence, $\B\spec\subseteq\A\spec$ implying
    \begin{equation*}
    \B\spec=\left\{ x\in\A\spec:\ \left\| x-z_{j}\right\| _{2}=r_{j}\text{  for all } j\right\}.
    \end{equation*}
    Since $c\left(r,Z\right){=}\projA{\spec}\left(z_{i}\right)$, for any $x{\in}\A\spec$, the vectors $\left(x{-}c\left(r,Z\right)\right)$ and $\left(z_{j}{-}c\spec\right)$ are orthogonal. Then, applying the Pythagorean Theorem to the equality condition in \eqref{eq:MBP_cons} gives
    \begin{equation*}
    r_{j}^{2} = \left\| x-z_i\right\|_{2}^{2} =\left\| x-c\left(r,Z\right)\right\| _{2}^{2}+\left\| z_{j}-c\spec\right\|_{2}^{2} \ \forall \ x\in\A\spec \text{ and } j=1,\dots,k.
    \end{equation*}
    We then deduce \eqref{eq:Bdef2} by subtracting $\left\| z_{j}{-}c\spec\right\|_{2}^{2}$ from the left and right sides, and applying the definition of $\rho\spec$ in \eqref{eq:margin}.
    
    \noindent \textit{Non-Emptiness of $\B\spec$:} By the first fact of Lemma \ref{lem:A}, the affine subspace $\A\spec$ is $(n{-}m{+}1)$-dimensional, so $\A\spec$ is a non-empty space. Thus, there exists $a\in\A\spec$ such that $a\neq c\spec$. Choose any such $a$ and consider the direction 
    \begin{equation*}
    d = \frac{1}{\left\| a-c\spec\right\|} \left(a-c\spec\right)
    \end{equation*}
    along with the vector
    \begin{equation*}
    v = c\spec+\rho\spec d.
    \end{equation*}
    Then, by construction,
    \begin{equation*}
    \left\| v-c\spec\right\| _{2}=\rho\spec.
    \end{equation*}
    Since $\A\spec$ contains both $c\spec$ and $d$, we deduce $v\in\A\spec$. Thus, $v$ satisfies all conditions in \eqref{eq:Bdef2}, so $v\in\B\spec$, and hence $\B\spec$ is non-empty.
    \end{proof}
    
    \begin{thm}\label{thm:MBP_sol}
    \textbf{(Local, Explicit Solution of $\boldsymbol{MBP\param}$)} For parameters $\param$ near some base parameters $\bparam$, the optimal solution to $\MBP\param$ is unique, varies smoothly with $\param$, and is given by\begin{equation}\label{eq:x}
    x\param = c\spec+\rho\spec u\param,
    \end{equation}
    where\begin{equation}\label{eq:u}
    u\param = \frac{1}{\left\| \projA{\spec}\left(y\right)-c\spec\right\| _{2}}\left(\projA{\spec}\left(y\right)-c\spec\right).
    \end{equation}
    In particular, $x\bparam = 0$.
    \end{thm}
    \begin{proof}
    \textit{The vector $u\param$ is well-defined:} To show the well-definedness of the expression in \eqref{eq:u}, observe the following at base parameters $\bparam$:
    
    \begin{align*}
    \projA{\bspec}\left(\by\right)-c\bspec &= \projA{\bspec}\left(-\sum_{j=1}^{\bk}\clam_{j}\bar{z}_{j}\right)-c\bspec &\text{ (By Definition \ref{def:base})}\\
    & = -\sum_{j=1}^{\bk}\clam_{j}\projA{\bspec}\left(\bar{z}_{j}\right)-c\bspec &\text{ (By linearity of $\A\bspec$ from Lemma \ref{lem:linA})}\\
    & = -\left(\sum_{j=1}^{\bk}\clam_{j}+1\right)c\bspec &\text{ (By definition of core \eqref{eq:core})}\\
    & \neq 0. &\text{ (By Definition \ref{def:base} and Lemma \ref{lem:linA})}.
    \end{align*}
    By the smoothness of the projection (third fact of Lemma \ref{lem:A}), we deduce $\projA{\spec}\left(y\right){-}c\spec{\neq}0$ for $\param$ sufficiently close to $\bparam$ as well. Thus, \eqref{eq:u} is well-defined. Moreover, the smoothness of the projection and \eqref{eq:u} further shows that $u\param$ is smooth with respect to $\param$.
    
    \noindent \textit{Derivation of $x\param$:} By Lemma \ref{lem:B}, $\B\spec\subseteq\A\spec$ so\begin{align*}
    x\param=\projB{\spec}\left(y\right) & =\projB{\spec}\left(\projA{\spec}\left(y\right)\right).
    \end{align*}
    Lemma \ref{lem:B} also states that $\B\spec$ is simply the $c\spec$-centered ball with radius $\rho\spec$ restricted to $\A\spec$. Thus, $\projB{\spec}\left(\projA{\spec}\left(y\right)\right)$ must be the (necessarily unique) point that
    \begin{enumerate}
        \item lies on the straight line from $c\spec$ to $\projA{\spec}\left(y\right)$, and
        \item is exactly the distance $\rho\spec$ from $c\spec$.
    \end{enumerate}
    Formalizing the above two criteria leads to the explicit-form
    \begin{equation*}
    x\param=c\spec+\frac{\rho\spec}{\left\| \projA{\spec}\left(y\right)-c\spec\right\| _{2}}\left(\projA{\spec}\left(y\right)-c\spec\right)=c\spec+\rho\spec u\param.
    \end{equation*}
    The smoothness of $x\param$ follows directly from the smoothness of the core (Lemma \ref{lem:core}), the margin (Lemma \ref{lem:margin}), and $u\param$ (shown above).
    
    \noindent \textit{Demonstration of $x\bparam{=}0$:} Observe that the Lagrangian of $\MBP\bparam$ is
    \[
    \L(x,\lambda) = \left\| x-\by\right\| _{2}^{2} + \sum_{j=1}^\bk \lambda_j \left\| x-\bz_{j}\right\| _{2}^{2},
    \]
    where $\lambda_j$ is the $j$-th element of a Lagrange multipliers vector $\lambda \in \R^\bk$. By Definitions \ref{def:MBP}-\ref{def:base}, the zero-vector is a feasible solution to $\MBP\bparam$. A feasible vector $x^*$ is the solution of $\MBP\bparam$ if and only if there exists $\lambda^*$ such that $x^*$ minimizes $\L(\cdot,\lambda^*)$. It is easy to check using Definition \ref{def:base} and first-order optimality conditions that $x^*{=}0$ is the unique minimizer of $\L(\cdot,\clam)$. Thus, $x\bparam{=}0$ is the unique solution to $\MBP\bparam$.
    \end{proof}
    
    \subsection{Lagrange Multipliers and Inequality Constraints}\label{sec:MBP_ineq}
    
    In this section, we will show that the optimal solutions of MBPs are unchanged when the equality constraints \eqref{eq:MBP_cons} are relaxed to inequalities. In particular, regardless of whether $\MBP\param$ is defined with equality or inequality constraints, note that the stationarity conditions remain the same. Therefore, there must exist some Lagrange multipliers, $\lambda\param\in\R^m$, such that\begin{equation}
    x\param-y+\sum_{j=1}^{\bk}\lambda_{j}\param\left[x\param-z_{j}\right]=0, \label{eq:MBP_FOC}
    \end{equation}
    where $\lambda_{j}\param$ denotes the $j$-th element of $\lambda\param$. The only difference is that the inequality version possesses additional dual feasibility conditions requiring the non-negativity of $\lambda\param$. Lemma \ref{lem:MBP_lam} will show that the $\lambda\param$ are unique and strictly positive.
    
    Then, for $\param$ sufficiently close to $\bparam$, Theorem \ref{thm:MBP_ineq} will show that the KKT conditions are sufficient to certify the optimality of $\left(x\param,\lambda\param\right)$ for inequality-constrained MBPs as well.
    
    \begin{lem}\label{lem:MBP_lam}\textbf{(Properties of MBP Lagrange Multipliers)} Suppose the parameters $\param$ are sufficiently close to base parameters $\bparam$ with corresponding base weights vector $\clam$. Consider the optimal solution $x\param$ of $\MBP\param$. Then, the following hold:
    \begin{enumerate}
        \item There exists a unique Lagrange multiplier vector, $\lambda\param\in\R_{++}^{\bk}$, satisfying the stationarity condition \eqref{eq:MBP_FOC}.
        \item The vector $\lambda\param$ varies smoothly with $\param$.
        \item At the base parameters, $\lambda\bparam=\clam$.
    \end{enumerate}
    \end{lem}
    \begin{proof}
    Since $x\param$ is optimal for $\MBP\param$, there must exist a Lagrange multipliers vector $\lambda\param\in\R^{\bk}$ satisfying \eqref{eq:MBP_FOC}. We first show $\lambda\param$ is unique. Finding Lagrange multipliers satisfying \eqref{eq:MBP_FOC} is equivalent to finding minimizers of the function
    \begin{equation*}
    F\left(\lambda\right)=\left\| x\param-y+X\lambda\right\| _{2}^{2},
    \end{equation*}
    where we define the matrix
    \[X=\begin{bmatrix}x\param{-}z_{1}, \dots, x\param{-}z_{\bk}\end{bmatrix}\in\R^{n\times\bk}.\]
    Then, a vector $\lambda$ minimizes $F$ if and only if it satisfies the below $\bk\times\bk$ linear system:
    \begin{equation*}
    X^{T}\left(x\param-y+X\lambda\right)=0\iff X^{T}X\lambda=X^{T}\left(y-x\param\right).
    \end{equation*}
    By Theorem \ref{thm:MBP_sol}, for $\param$ close to $\bparam$, the MBP-solution $x\param$ is close to zero. Therefore, columns of $X$ are close to $\left\{-\bz_{j}\right\}_{j=1}^{\bk}$, which are linearly independent (Definition \ref{def:base}). This implies $X^{T}X$ is invertible, and the minimizer of $F$ is unique and explicitly given by 
    \begin{equation*}
    \lambda\param = \argmin_\lambda F\left(\lambda\right)=\left(X^{T}X\right)^{-1}X^{T}\left(y-x\param\right).
    \end{equation*}
    The smoothness of $\lambda\param$ with respect to $\param$ follows from the smoothness of $x\param$ (Theorem \ref{thm:MBP_sol}) and the smoothness of the above matrix operations. 
    
    Finally, it follows directly from Definition \ref{def:base} and Theorem \ref{thm:MBP_sol} that
    \begin{equation*}
    x\bparam-\by+\sum_{j=1}^{\bk}\clam_{j}\left[x\bparam-\bz_{j}\right]=0,
    \end{equation*}
    which implies the Lagrange multiplier vector at the base parameters is uniquely $\lambda\bparam=\clam$.
    \end{proof}
    
    \begin{thm}\textbf{(Inequality-Constrained MBP)}\label{thm:MBP_ineq}
    Consider the following inequality-constrained MBP:
    \begin{align*}
    \min_{x\in\R^{n}}\ & \left\| x-y\right\| _{2}^{2}\\
    \mathrm{s.t.}\ & \left\| x-z_{j}\right\| _{2}^{2} \leq r_{j}^{2}\ \forall\ j=1,\dots,\bk,
    \end{align*}
    and denote it as $\MBPl\param$. For parameters $\param$ sufficiently close to $\bparam$, $\MBPl\param$ is feasible and possesses the same unique optimal solution, $x\param$, and unique Lagrange multipliers, $\lambda\param$, as $\MBP\param$. Consequently, all constraints of $\MBPl\param$ are tight. Moreover, the optimal solution and multipliers satisfy\begin{equation}
    x\param = \mathrm{O}\left(\|\param-\bparam\|\right)\text{ and } \lambda\param = \clam + \mathrm{O}\left(\|\param-\bparam\|\right).\label{eq:MBP_Onotation}
    \end{equation}
    \end{thm}
    \begin{proof}
    The feasibility of $\MBPl\param$ follows directly from Lemma \ref{lem:B} since $\B\param$ is a subset of the feasible region of $\MBPl\param$.
    
    $\MBPl\param$ possesses a unique minimizer since it minimizes a strongly convex objective on a convex feasible region. Moreover, checking the KKT conditions is sufficient for demonstrating optimality because the $\|x - z_j\|_2^2$-term in the constraints is convex and differentiable for all $j$. Therefore---denoting the unique pair of optimal solutions and Lagrange multipliers of $\MBP\param$ as $\left(x\param, \lambda\param\right)$---it suffices to check that $\left(x\param, \lambda\param\right)$ satisfies the KKT conditions of $\MBPl\param$:
    \begin{enumerate}
        \item \textbf{(Stationarity)} Follows from $\MBPl\param$ and $\MBP\param$ having the same stationarity condition \eqref{eq:MBP_FOC}.
        \item \textbf{(Primal Feasibility)} Follows from $\B\param$ being a subset of the feasible region of $\MBPl\param$.
        \item \textbf{(Dual Feasibility)} The inequality constraints of $\MBPl\param$ require the non-negativity of all Lagrange multipliers. By Definition \ref{def:base}, the elements of $\clam$ are strictly-positive. Thus, by Lemma \ref{lem:MBP_lam}, the elements of $\lambda\param$ are also strictly-positive for $\param$ sufficiently close to $\bparam$.
        \item \textbf{(Complementary Slackness)} \sloppy By construction, $x\param$ satisfies all the inequality constraints of $\MBPl\param$ with equality. Correspondingly, $\lambda\param$ are all strictly positive.
    \end{enumerate}
    Therefore, $\left(x\param, \lambda\param\right)$ is the unique pair of optimal solutions and Lagrange multipliers for $\MBPl\param$ as well. The characterizations in \eqref{eq:MBP_Onotation} directly follows from Theorem \ref{thm:MBP_sol} and Lemma \ref{lem:MBP_lam} as well as recalling that smooth functions are locally Lipschitz.
    \end{proof}
    
    \subsection{Survey Descent as a Inequality Constrained MBP}\label{sec:MBP_SD}
    
    We can reformulate Survey Descent subproblems into inequality-constrained MBPs. In particular, for all $i$, we can rewrite the $i$-th Survey Descent subproblem in Observation \ref{obs:SD_MOS} as
    
    \begin{gather}\begin{aligned}\min_{x}\ & \left\| (x - s_{i}) - y_{i}\left(\S\right) \right\| _{2}^{2}\\
    \mathrm{s.t.}\ & \left\| (x-s_{i})-z_{ij}\left(\S\right)\right\| _{2}^{2}\leq r_{ij}^{2}\left(\S\right)\ \forall\ j\neq i,
    \end{aligned}
    \label{eq:proj_form}
    \end{gather}
    where, for all $i$ and $j$, we define the following mappings using the components $\{f_i\}_{i=1}^k$ from Definition \ref{def:strong_max}:
    
    \begin{gather}\begin{aligned}
    y_{i}\left(\S\right)& \equiv-\frac{1}{L}\nabla f_{i}\left(s_{i}\right),\\
    z_{ij}\left(\S\right) & \equiv s_j - s_i  + y_{j}\left(\S\right)-y_{i}\left(\S\right),\\
    r_{ij}^{2}\left(\S\right) & \equiv\frac{2}{L}\left[f_{i}\left(s_{i}\right)-f_{j}\left(s_{j}\right)+\nabla f_{i}\left(s_{i}\right)^{T}\left(s_{j}-s_{i}\right)\right]+\left\|y_{j}\left(\S\right)-y_{i}\left(\S\right)\right\| _{2}^{2}.
    \end{aligned}\label{eq:MBP_mappings}
    \end{gather}
    \begin{observe}\textbf{(Survey Descent as Inequality-Constrained MBP)} \label{obs:SD_as_MBP} Consider any $k$-point survey $\S$ and any $i{=}1,\dots,k$. After the changes-of-variables $x \leftarrow x{-} s_i$, the minimization problem in \eqref{eq:proj_form} is an inequality-constrained MBP  with $\m = k{-}1$ and parameters 
    
    $$\param = \left(y_{i}\left(\S\right), \left[ z_{ij}\left(\S\right)\right]_{j\neq i}, \left[r_{ij}\left(\S\right)\right]_{j\neq i},\left[ s_{j}\right]_{j\neq i}\right),$$
    where, for any $j$-indexed scalars $\left\{v_j\right\}_{j=1}^{k}$, we use the notation $\left[v_j\right]_{j \neq i}{\equiv} \left[v_1,\dots,v_{i-1},v_{i+1},\dots,v_{k}\right]^T$ to denote the length $k{-}1$ vector excluding $v_i$.
    Furthermore, the mappings \eqref{eq:MBP_mappings}  evaluated at the base survey\footnote{Evaluating of the mappings \eqref{eq:MBP_mappings} at $\bX$ is a technical step. $\bX$ is not a valid survey (Definition \ref{def:validS}), so it does not correspond to an actual Survey Descent iteration.} $\bX$ (Definition \ref{def:validS}) form a set of base parameters
    \begin{equation}
        \bparam = \left(y_{i}\left(\bX\right), \left[ z_{ij}\left(\bX\right)\right]_{j\neq i}, \left[r_{ij}\left(\bX\right)\right]_{j\neq i}, [\bx]_{j\neq i}\right),\label{eq:MBP_SD_base}
    \end{equation}
    with corresponding weights $\{\blam_j\}_{j\neq i}$ from \eqref{eq:lam_sum1}-\eqref{eq:lam_sum0}, and where $[\bx]_{j\neq i}$ is the $k{-}1$ length vector with all $\bx$ entries.
    \end{observe}
    We then deduce our desired theorem.
    \PBmbp*
    \begin{proof}
    It is easy to check from the mapping definitions \eqref{eq:MBP_mappings} and the Taylor expansion \eqref{eq:taylor} that 
    
    $$\left\|\left(y_{i}\left(\S\right), \left[ z_{ij}\left(\S\right)\right]_{j\neq i}, \left[r_{ij}\left(\S\right)\right]_{j\neq i},\left[ s_{j}\right]_{j\neq i}\right) - \left(y_{i}\left(\bX\right), \left[ z_{ij}\left(\bX\right)\right]_{j\neq i}, \left[r_{ij}\left(\bX\right)\right]_{j\neq i}, [\bx]_{j\neq i}\right)\right\| = \mathrm{O}\left(\|\S-\bX\|\right).$$
    The result then follows directly from Observation \ref{obs:SD_MOS}, Observation \ref{obs:SD_as_MBP}, and Theorem \ref{thm:MBP_ineq}.
    \end{proof}
    
    Observe that the proofs of Lemmas \ref{lem:A}-\ref{lem:MBP_lam} explicitly construct both $x\param$ and $\lambda\param$ using simple linear algebra operations and a scalar square-root. Thus, combining these lemmas with Observation \ref{obs:SD_as_MBP}, we deduce the below computationally-efficient procedure for solving Survey Descent subproblems within the setting of Theorem \ref{thm:SD_MBP}. Moreover, recall that $\nabla f(s_i){=}\nabla f_i(s_i)$ since $\S$ is valid, so the procedure requires only a first-order oracle on $f$ without need to access its $\{f_i\}_{i=1}^k$ components.
    \begin{observe}
    \textbf{(Survey Descent with Linear Algebra)} Assume Setting \ref{setting:localSD}. Consider $i{=}k$ without loss of generality. We can compute the solution, $s_k^+(\S)$, and Lagrange multipliers, $\{\lambda_j(\S)\}_{j=1}^{k-1}$, of the Survey Descent subproblem $(P_k^\S)$ using the following procedure:
    \begin{enumerate}
    \item Compute
    \begin{align*}
        y_j& \leftarrow -\frac{1}{L}\nabla f\left(s_{j}\right) & \forall \ j = 1,\dots,k,\\
        z_j & \leftarrow s_j - s_k  + y_{j}-y_{k} & \forall \ j=1,\dots,k{-}1,\\
        r_{j}^{2} & \leftarrow \frac{2}{L}\left[f\left(s_{k}\right)-f\left(s_{j}\right)+\nabla f\left(s_{k}\right)^{T}\left(s_{j}-s_{k}\right)\right]+\left\|y_{j}-y_{k}\right\| _{2}^{2} & \forall \ j=1,\dots,k{-}1,\\
        M &\leftarrow 2\begin{bmatrix}\left(z_{2}-z_{1}\right)^{T}\\
        \vdots\\
        \left(z_{k-1}-z_{1}\right)^{T}
        \end{bmatrix}, \\
        b &\leftarrow\begin{bmatrix}\left(\left\| z_{2}\right\| _{2}^{2}-r_{2}^{2}\right)-\left\| z_{1}\right\| _{2}^{2}+r_{1}^{2}\\
        \vdots\\
        \left(\left\| z_{k-1}\right\| _{2}^{2}-r_{k-1}^{2}\right)-\left\| z_{1}\right\| _{2}^{2}+r_{1}^{2}
        \end{bmatrix}.
    \end{align*}
    \item Define the projection function
    \begin{equation*}
    \P\left(w\right)\equiv\left(I-M\left(M^{T}M\right)^{-1}M^{T}\right)w+M\left(M^{T}M\right)^{-1}b.
    \end{equation*}
    We would implement the above projection with a QR-factorization instead of explicitly taking the matrix inverse, which can be numerically unstable (see, for example, \citet[Chapter 14.1]{higham2002accuracy}).
    \item Choose and fix any $\bj$ from $1,\dots,k{-}1$.
    \item Compute 
    \begin{align*}
        c &\leftarrow\P\left(z_{\bj}\right),\\
        \rho &\leftarrow\sqrt{r_{\bj}^{2}-\left\|z_{\bj}-c\right\| _{2}^{2}},\\
        x &\leftarrow c+\frac{\rho}{\left\|\P\left(y_k\right)-c\right\| _{2}}\left(\P\left(y_k\right)-c\right),\\
        X &\leftarrow \begin{bmatrix}x-z_{1}, & \dots, & x-z_{k-1}\end{bmatrix},\\
        \lambda^{\prime} &\leftarrow \left(X^{T}X\right)^{-1}X^{T}\left(y_k-x\right).
    \end{align*}
    Similar to step 2, we would compute $\lambda^{\prime}$ using a QR-factorization for numerical stability.
    \item Output
    \begin{align*}
        s_k^+(\S) &\leftarrow x + s_k, \\
        \lambda_j(\S) &\leftarrow \lambda_j^{\prime} \qquad \forall \ j=1,\dots,k-1,
    \end{align*}
    where $\lambda_j^{\prime}$ denotes the $j$-th element of the vector $\lambda^\prime$.
    \end{enumerate}
    The general procedure for solving $(P_i^\S)$, for any $i=1,\dots,k$, follows from re-indexing.
    \end{observe}
    
    The variable $x$ in steps 4 and 5 is directly related to Observation \ref{obs:SD_as_MBP}, which formulates the $i$-th Survey Descent subproblem as an inequality-constrained MBP after the variable-shift, $x{\gets}x{-}s_i$. In step 4, we compute the solution, $x$, of the inequality-constrained MBP using the shifted-variables. Step 5 then reverses the variable shift to obtained the survey point update.
\fi 

\end{document}